\documentclass{TTP_DSL2006}

\usepackage{graphicx,graphics}
\usepackage[intlimits]{amsmath}
\usepackage{amssymb}
\usepackage{exscale}
\usepackage{times}
\usepackage{subfigure}
\usepackage{multirow}
\usepackage[abs]{overpic}

\renewcommand{\large}{\fontsize{14}{18pt}\selectfont}
\renewcommand{\small}{\fontsize{11}{13.6pt}\selectfont}

\newcommand{\titleformat}{\sffamily\bfseries \large}						
\newcommand{\authorformat}{\sffamily \large}							
\newcommand{\keywordsformat}{\noindent \small \sffamily}				
\newcommand{\abstractformat}{\noindent \textbf}						
\newcommand{\contentformat}{\rmfamily \normalsize\vspace{18pt}}			
\newcommand{\email}{\sffamily \small \vspace{-8pt}}						
\renewcommand{\subsection}{\textbf}	

\newcommand{\Eref}[1]{Equation (\ref{#1})}
\newcommand{\fref}[1]{Figure (\ref{#1})}
\newcommand{\Erefs}[1]{Equations (\ref{#1})}

\newcommand{\xx}{\mathbf{x}}

\newcommand{\bm}{\mathbf{M}}
\newcommand{\bn}{\mathbf{N}}
\newcommand{\DD}{\mathbf{D_b}}

\newcommand{\bveps}{\boldsymbol{\varepsilon}}
\newcommand{\uu}{\mathbf{u}}
\newcommand{\BB}{\mathbf{B}}
\newcommand{\qq}{\mathbf{q}}
\newcommand{\aaa}{\mathbf{a}}
\newcommand{\bb}{\mathbf{b}}

\newcommand{\cc}{\mathbf{c}}
\newcommand{\rmd}{\rm{d}}

\begin{document}

\title{\titleformat Vibration of functionally graded material plates with cutouts \& cracks in thermal environment}

\author{\authorformat Ahmad Akbari Rahimabadi\inst{1}$^{,\rm{a}}$\text{,} Sundararajan Natarajan \inst{2}$^{,\rm{b}}$\text{,} and St\'ephane PA Bordas \inst{1}$^{,\rm{c}}$}


\institute{\sffamily Institute of Mechanics and Advanced Materials, Cardiff School of Engineering, Cardiff University, Wales, CF24 3AA, UK. \and School of Civil \& Environmental Engineering, The University of New South Wales, Sydney, NSW 2052, Australia. }

\maketitle

\begin{center}
\email{ $^{\rm a}$ahmad.akbari.r@gmail.com, $^{\rm b}$sundararajan.natarajan@gmail.com, $^{\rm c}$stephane.bordas@gmail.com}
\end{center}

\keywordsformat{{\textbf{Keywords:}}  Vibration, cutouts, cracks, Reissner-Mindlin plate, extended finite element method.}

\contentformat

\abstractformat{Abstract.} In this paper, the effect of a centrally located cutout (circular and elliptical) and cracks emanating from the cutout on the free flexural vibration behaviour of functionally graded material plates in thermal environment is studied. The discontinuity surface is represented independent of the mesh by exploiting the partition of unity method framework. A Heaviside function is used to capture the jump in the displacement across the discontinuity surface and asymptotic branch functions are used to capture the singularity around the crack tip. An enriched shear flexible 4-noded quadrilateral element is used for the spatial discretization. The properties are assumed to vary only in the thickness direction. The effective properties of the functionally graded material are estimated using the Mori-Tanaka homogenization scheme and the plate kinematics is based on the first order shear deformation theory. The influence of the plate geometry, the geometry of the cutout, the crack length, the thermal gradient and the boundary conditions on the free flexural vibration is numerically studied. 

\section{Introduction}
\vspace{-6pt}

In recent years, a new class of engineered material, the functionally graded materials (FGMs) has attracted researchers to investigate its structural behaviour. The FGMs are microscopically inhomogeneous and the mechanical and the thermal properties vary \textit{smoothly and continuously} from one surface to another. FGMs combine the best properties of its constituents. Typically, the FGMs are made from a mixture of ceramic and metal. The ceramic constituent provides thermal stability due to its low thermal conductivity, whilst the metallic constituent provides structural stability. FGMs eliminate the sharp interfaces existing in laminated composites with a gradient interface and are considered to be an alternative in many engineering sectors such as the aerospace industry, biomechanics industry, nuclear industry, tribology, optoelectronics and other high performance applications where the structural member is exposed to high thermal gradient in addition to mechanical load. 

The static and the dynamic characteristics have been studied in detail by many researchers using different plate theories, for example, first order shear deformation theory (FSDT)~\cite{Reddy2000,Yang2002,Sundararajan2005}, second and other higher order accurate theory~\cite{Qian2004a,Ferreira2006,natarajanmanickam2012} have been used to describe the plate kinematics. Existing approaches in the literature to study plate and shell structures made up of FGMs uses finite element method (FEM) based on Lagrange basis functions~\cite{Reddy2000,ganapathiprakash2006,Sundararajan2005}, meshfree methods~\cite{Qian2004a,Ferreira2006} and recently Valizadeh \textit{et al.,}~\cite{valizadehnatarajan2013} used non-uniform rational B-splines based FEM to study the static and the dynamic characteristics of FGM plates in thermal environment. Akbari \textit{et al.,}~\cite{R2010} studied two-dimensional wave propagation in functionally graded solids using the meshless local Petrov-Galerkin method. Huang \textit{et al.,}~\cite{Huang2011} proposed solutions for the free vibration of side-cracked FGM thick plates based on Reddy's third-order shear deformation theory using Ritz technique. Kitipornchai~\textit{et al.,}~\cite{Kitipornchai2009} studied nonlinear vibration of edge cracked functionally graded Timoshenko beams using Ritz method. Yang~\textit{et al.,}~\cite{Yang2010} studied the nonlinear dynamic response of a functionally graded plate with a through-width crack based on Reddy's third-order shear deformation theory using a Galerkin method. Natarajan \textit{et al.,}~\cite{natarajanbaiz2011,natarajanbaiz2011a} and Baiz \textit{et al.,}~\cite{baiznatarajan2011} studied the influence of the crack length on the free flexural vibrations of FGM plates using the XFEM and smoothed XFEM, respectively. The above list is no way comprehensive and interested readers are referred to the literature and references therein and a recent review paper by Jha and Kant~\cite{jhakant2013} on FGM plates. 

Plates with cutouts are extensively used in transport vehicle structures. Cutouts are made to lighten the structure, for ventilation, to provide accessibility to other parts of the structures and for altering the resonant frequency. Therefore, the natural frequencies of plates with cutouts are of considerable interest to designers of such structures. Most of the earlier investigations on plates with cutouts have been confined to isotropic plates~\cite{paramasivam1973,aliatwal1980,huangsakiyama1999} and laminated composites~\cite{reddy1982,sivakumariyengar1998}. Recently, Janghorban and Zare~\cite{janghorbanzare2011} studied the influence of cutout on the fundamental frequency of FGM plates in thermal environment using FEM. Their study was restricted to a limited number of configurations, because the mesh has to conform to the geometry. In this study, we present a framework that provides flexibility to handle internal discontinuities.

In this paper, we study the influence of a centrally located cutout and cracks emanating from the cutouts on the natural frequencies of FGM plates in thermal environment. Circular and elliptical cutouts are considered for the study. A structured quadrilateral mesh is used and the cutouts are modelled independent of the mesh within the extended finite element (XFEM) framework. A systematic parametric study is carried out to bring the effect of the gradient index $(k)$, the thermal gradient $(\Delta T = T_c - T_m)$, the geometry of the cutout on the free flexural vibration behaviour of FGM plates with different boundary conditions. \\


\section{Theoretical Background}
\label{theory}
\begin{figure}[htpb]
\centering
\includegraphics[scale=0.6]{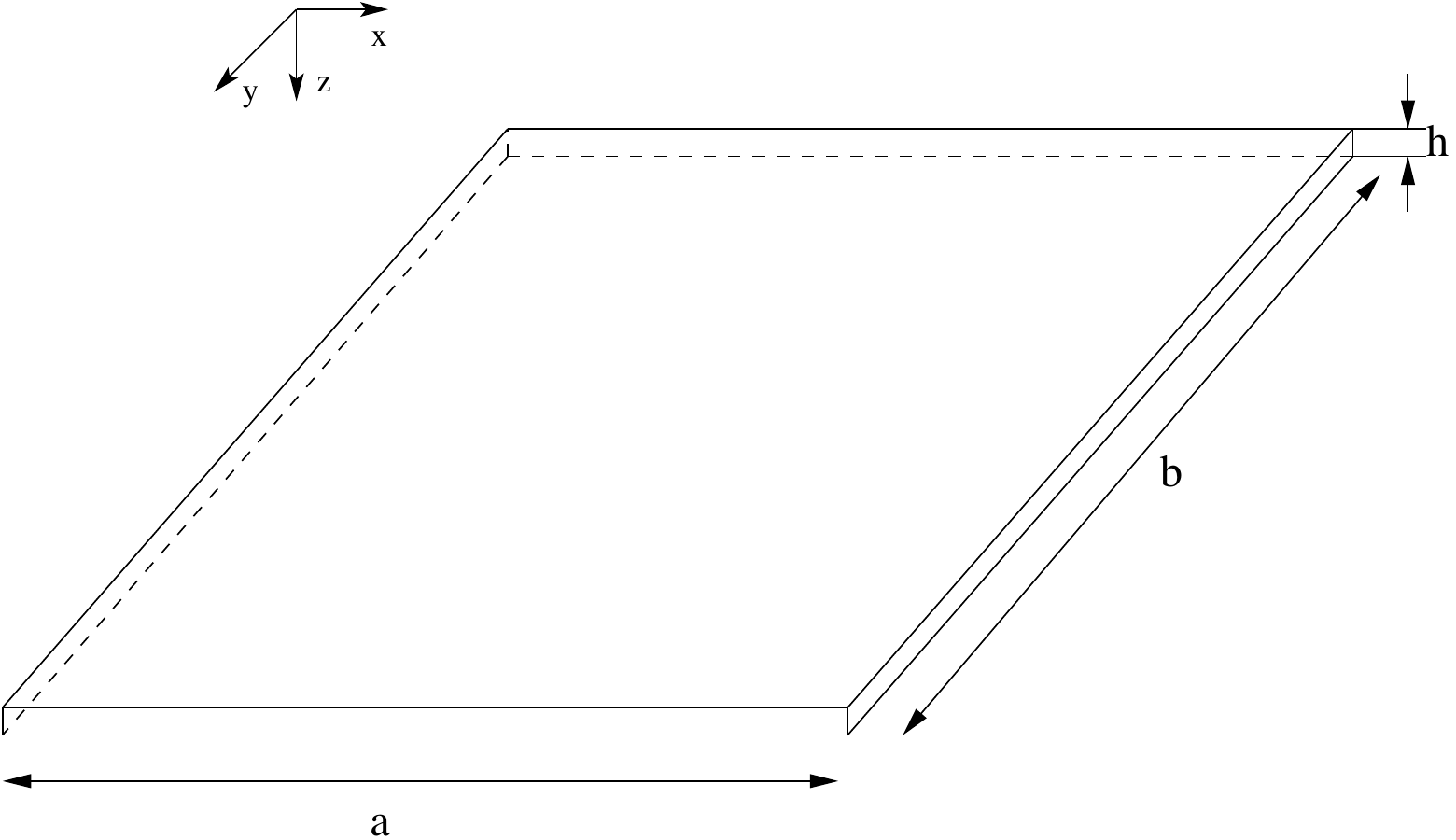}
\caption{Coordinate system of a rectangular plate.}
\label{fig:platefig}
\end{figure}

\textbf{Reissner-Mindlin plate theory} The Reissner-Mindlin plate theory, also known as the first order shear deformation plate theory (FSDT) is an extension of the classical plate theory (or, the Kirchhoff - Love plate theory). The FSDT takes into account the shear deformations through the thickness and are intended for thick plates in which the normal to the medium surface remains straight but not necessarily perpendicular to the medium surface. Using the Mindlin formulation, the displacements $u,v,w$ at a point $(x,y,z)$ in the plate (see \fref{fig:platefig}) from the medium surface are expressed as functions of the mid-plane displacements $u_o,v_o,w_o$ and independent rotations $\theta_x,\theta_y$ of the normal in $yz$ and $xz$ planes, respectively, as

\begin{align}
u(x,y,z,t) &= u_o(x,y,t) + z \theta_x(x,y,t) \nonumber \\
v(x,y,z,t) &= v_o(x,y,t) + z \theta_y(x,y,t) \nonumber \\
w(x,y,z,t) &= w_o(x,y,t) 
\label{eqn:displacements}
\end{align}
where $t$ is the time. The strains in terms of mid-plane deformation can be written as:
\begin{equation}
\bveps  = \left\{ \begin{array}{c} \bveps_p \\ \mathbf{0} \end{array} \right \}  + \left\{ \begin{array}{c} z \bveps_b \\ \bveps_s \end{array} \right\} 
\label{eqn:strain1}
\end{equation}
The midplane strains $\bveps_p$, the bending strains $\bveps_b$ and the shear strains $\varepsilon_s$ in \Eref{eqn:strain1} are written as:

\begin{equation}
\renewcommand{\arraystretch}{1.5}
\bveps_p = \left\{ \begin{array}{c} u_{o,x} \\ v_{o,y} \\ u_{o,y}+v_{o,x} \end{array} \right\}, \hspace{1cm}
\renewcommand{\arraystretch}{1.5}
\bveps_b = \left\{ \begin{array}{c} \theta_{x,x} \\ \theta_{y,y} \\ \theta_{x,y}+\theta_{y,x} \end{array} \right\} , \hspace{1cm}
\renewcommand{\arraystretch}{1.5}
\bveps_s = \left\{ \begin{array}{c} \theta _x + w_{o,x} \\ \theta _y + w_{o,y} \end{array} \right\}.
\renewcommand{\arraystretch}{1.5}
\end{equation}
where the subscript `comma' represents the partial derivative with respect to the spatial coordinate succeeding it. The membrane stress resultants $\bn$ and the bending stress resultants $\bm$ can be related to the membrane strains, $\bveps_p$ and the bending strains $\bveps_b$ through the following constitutive relations:

\begin{equation}
\bn = \left\{ \begin{array}{c} N_{xx} \\ N_{yy} \\ N_{xy} \end{array} \right\} = \mathbf{A} \bveps_p + \BB \bveps_b - \bn^\textup{th}, \hspace{0.5cm}
\bm = \left\{ \begin{array}{c} M_{xx} \\ M_{yy} \\ M_{xy} \end{array} \right\} = \BB \bveps_p + \DD \bveps_b - \bm^\textup{th}.
\end{equation}
where the matrices $\mathbf{A} = A_{ij}, \BB= B_{ij}$ and $\DD = D_{ij}; (i,j=1,2,6)$ are the extensional, the bending-extensional coupling and the bending stiffness coefficients and are defined as:

\begin{equation}
\left\{ A_{ij}, ~B_{ij}, ~ D_{ij} \right\} = \int_{-h/2}^{h/2} \overline{Q}_{ij} \left\{1,~z,~z^2 \right\}~\rmd z
\end{equation}
The thermal stress resultants, $\bn^{\textup{th}}$ and the moment resultants $\bm^{\textup{th}}$ are:

\begin{equation}
\left\{ \bn^{\rm th}, \bm^{\rm th} \right\} = \left\{ \begin{array}{c} N_{xx}^{\textup{th}}, M_{xx}^{\textup{th}} \\ N_{yy}^{\textup{th}}, M_{yy}^{\textup{th}}  \\ N_{xy}^{\textup{th}}, M_{xy}^{\textup{th}}  \end{array} \right\} = \int_{-h/2}^{h/2} \overline{Q}_{ij} \alpha(z,T) \left\{ \begin{array}{c} 1 \\ 1 \\ 0 \end{array} \right\} ~ \Delta T(z)  \left\{1,z \right\} ~ \rmd z
\end{equation}
where the thermal coefficient of expansion $\alpha(z,T)$ is given by \Eref{eqn:thermalcondalpha} and $\Delta T(z) = T(z)-T_o$ is the temperature rise from the reference temperature $T_o$ at which there are no thermal strains. Similarly, the transverse shear force $\mathbf{Q} = \{Q_{xz},Q_{yz}\}$ is related to the transverse shear strains $\bveps_s$ through the following equation:

\begin{equation}
\mathbf{Q} = \mathbf{E} \bveps_s
\end{equation}
 where $\mathbf{E} = E_{ij} = \int_{-h/2}^{h/2} \overline{Q}_{ij} \upsilon_i \upsilon_j~dz;~ (i,j=4,5)$ is the transverse shear stiffness coefficient, $\upsilon_i, \upsilon_j$ is the transverse shear coefficient for non-uniform shear strain distribution through the plate thickness. The stiffness coefficients $\overline{Q}_{ij}$ are defined as:

\begin{eqnarray}
\overline{Q}_{11} = \overline{Q}_{22} = {E(z,T) \over 1-\nu^2}; \hspace{1cm} \overline{Q}_{12} = {\nu E(z,T) \over 1-\nu^2}; \hspace{1cm} \overline{Q}_{16} = \overline{Q}_{26} = 0; \nonumber \\
\overline{Q}_{44} = \overline{Q}_{55} = \overline{Q}_{66} = {E(z,T) \over 2(1+\nu) }; \hspace{1cm} \overline{Q}_{45} = \overline{Q}_{54} = 0.
\end{eqnarray}

\noindent where the modulus of elasticity $E(z,T)$ and Poisson's ratio $\nu$ are given by \Eref{eqn:young}. The strain energy function $U$ is given by:

\begin{equation}
U(\boldsymbol{\delta}) = {1 \over 2} \int_{\Omega} \left\{ \bveps_p^{\textup{T}} \mathbf{A} \bveps_p + \bveps_p^{\textup{T}} \mathbf{B} \bveps_b + 
\bveps_b^{\textup{T}} \mathbf{B} \bveps_p + \bveps_b^{\textup{T}} \mathbf{D} \bveps_b +  \bveps_s^{\textup{T}} \mathbf{E} \bveps_s -  \bveps_p^{\textup{T}} \bn ^{\rm th}- \bveps_b^{\textup{T}} \bm^{\rm th} \right\}~ \rmd\Omega
\label{eqn:potential}
\end{equation}

\noindent where $\boldsymbol{\delta} = \{u,v,w,\theta_x,\theta_y\}$ is the vector of the degrees of freedom associated to the displacement field in a finite element discretization. Following the procedure given in~\cite{Rajasekaran1973}, the strain energy function $U$ given in~\Eref{eqn:potential} can be rewritten as:

\begin{equation}
U(\boldsymbol{\delta}) = {1 \over 2}  \boldsymbol{\delta}^{\textup{T}} \mathbf{K}  \boldsymbol{\delta}
\label{eqn:poten}
\end{equation}
where $\mathbf{K}$ is the linear stiffness matrix. The kinetic energy of the plate is given by:

\begin{equation}
T(\boldsymbol{\delta}) = {1 \over 2} \int_{\Omega} \left\{p (\dot{u}_o^2 + \dot{v}_o^2 + \dot{w}_o^2) + I(\dot{\theta}_x^2 + \dot{\theta}_y^2) \right\}~\rmd \Omega
\label{eqn:kinetic}
\end{equation}
where $p = \int_{-h/2}^{h/2} \rho(z)~\rmd z$, $I = \int_{-h/2}^{h/2} z^2 \rho(z)~\rmd z$ and $\rho(z)$ is the mass density that varies through the thickness of the plate. The external work due to the in-plane stress resultants $(N_{xx}^{\textup{th}}, N_{yy}^{\textup{th}}, N_{xy}^{\textup{th}})$ developed in the plate under the thermal load is:

\begin{equation}
\begin{split}
V(\boldsymbol{\delta}) = \int_{\Omega} \left\{ {1 \over 2} \left[ N_{xx}^{\textup{th}} w_{o,x}^2 + N_{yy}^{\textup{th}} w_{o,y}^2 +
2 N_{xy}^{\textup{th}} w_{o,x} w_{o,y}\right] \right. \\
\left. + {h^2 \over 24} \left[ N_{xx}^{\textup{th}} \left( \theta_{x,x}^2 + \theta_{y,x}^2 \right) + N_{yy}^{\textup{th}} \left( \theta_{x,y}^2 + \theta_{y,y}^2 \right) + 2N_{xy}^{\textup{th}} \left( \theta_{x,x} \theta_{x,y} + \theta_{y,x} \theta_{y,y} \right) \right] \right\} ~\rmd \Omega
\end{split}
\label{eqn:potthermal}
\end{equation}
Upon substituting \Erefs{eqn:poten} - (\ref{eqn:potthermal}) in Lagrange's equations of motion, the following governing equation is obtained:

\begin{equation}
\mathbf{M} \ddot{\boldsymbol{\delta}} + (\mathbf{K}+\mathbf{K}_{\rm G}) \boldsymbol{\delta} = \mathbf{0}
\label{eqn:govereqn}
\end{equation}
where $\mathbf{M}$ is the consistent mass matrix. After substituting the characteristic of the time function~\cite{Ganapathi1991} $\ddot{\boldsymbol{\delta}} = -\omega^2 \boldsymbol{\delta}$, the following algebraic equation is obtained:

\begin{equation}
\left( (\mathbf{K}+\mathbf{K}_{\rm G})  - \omega^2 \mathbf{M}\right) \boldsymbol{\delta} = \mathbf{0}
\label{eqn:finaldiscre}
\end{equation}
where $\mathbf{K}$ and $\mathbf{K}_{\rm G}$ are the stiffness matrix and the geometric stiffness matrix due to thermal loads, respectively, $\omega$ is the natural frequency. The plate is subjected to a temperature field. So, the first step in the solution process is to compute the in-plane stress resultants due to the temperature field. These will then be used to compute the stiffness matrix of the system and then the frequencies are computed for the system. \\

\textbf{Functionally Graded Material} A functionally graded material (FGM) rectangular plate (length $a$, width $b$ and thickness $h$), made by mixing two distinct material phases: a metal and a ceramic is considered with coordinates $x,y$ along the in-plane directions and $z$ along the thickness direction (see \fref{fig:platefig}). The material on the top surface $(z=+h/2)$ of the plate is ceramic and is graded to metal at the bottom surface of the plate $(z=-h/2)$ by a power law distribution. The homogenized material properties are computed using the Mori-Tanaka Scheme~\cite{Mori1973,Benvensite1987}. 

\textbf{Estimation of mechanical and thermal properties}
Based on the Mori-Tanaka homogenization method, the effective bulk modulus $K$ and shear modulus $G$ of the FGM are evaluated as~\cite{Mori1973,Benvensite1987,Cheng2000,Qian2004}

\begin{equation}
{K - K_m \over K_c - K_m} = {V_c \over 1+(1-V_c){3(K_c - K_m) \over 3K_m + 4G_m}}, \hspace{0.4cm}
{G - G_m \over G_c - G_m} = {V_c \over 1+(1-V_c){(G_c - G_m) \over G_m + f_1}},
\label{eqn:bulkshearmodulus}
\end{equation}
where
\begin{equation}
f_1 = {G_m (9K_m + 8G_m) \over 6(K_m + 2G_m)}.
\end{equation}

\begin{figure}[htpb]
\centering
\includegraphics[scale=0.65]{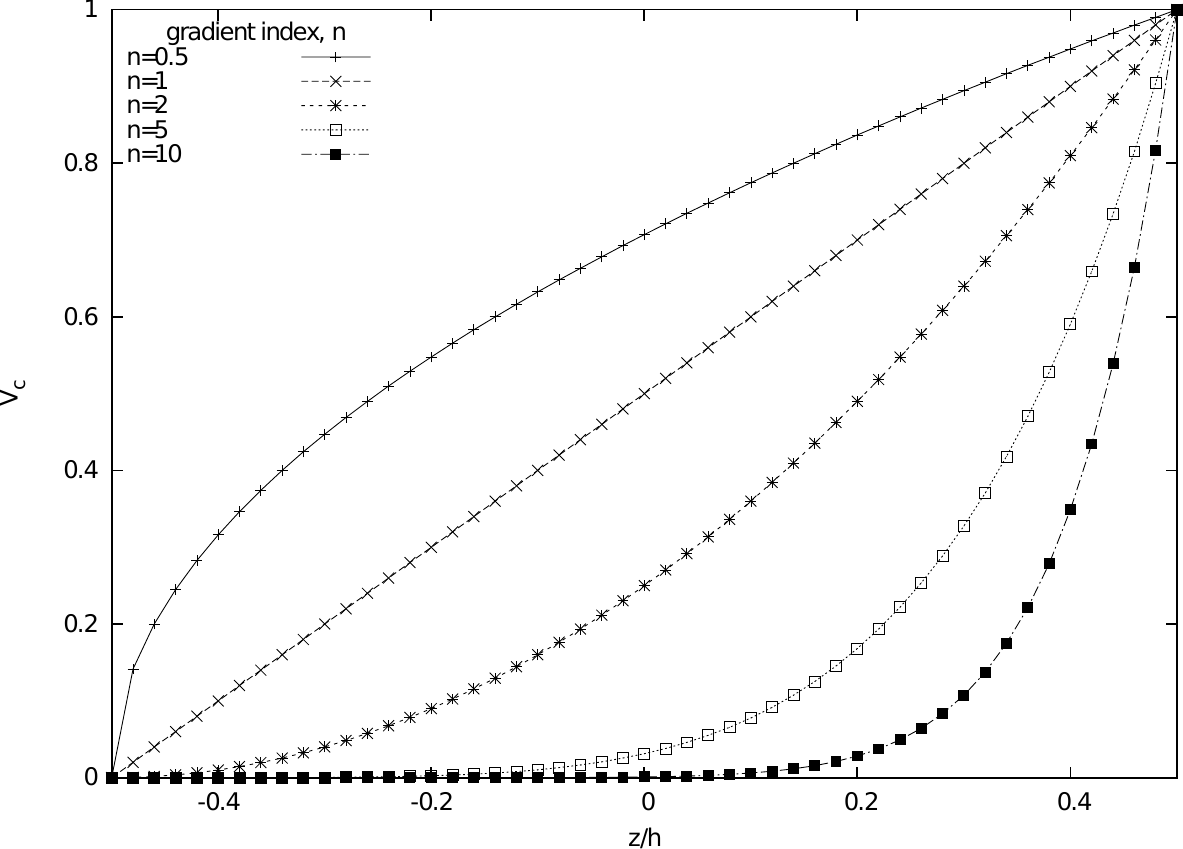}
\caption{Through thickness variation of volume fraction}
\label{fig:volfrac}
\end{figure}

Here, $V_i~(i=c,m)$ is the volume fraction of the phase material. The subscripts $c$ and $m$ refer to the ceramic and the metal phases, respectively. The volume fractions of the ceramic and the metal phases are related by $V_c + V_m = 1$ and $V_c$ is expressed as
\begin{equation}
V_c(z) = \left( {2z + h \over 2h} \right)^k, \hspace{0.2cm}  k \ge 0
\label{eqn:volFrac}
\end{equation}
where $k$ in \Eref{eqn:volFrac} is the volume fraction exponent, also referred to as the gradient index. \fref{fig:volfrac} shows the variation of the volume fraction of ceramic constituent in the thickness direction $z$ for a FGM plate. The effective Young's modulus $E$ and Poisson's ratio $\nu$ can be computed from the following expressions:
\begin{equation}
E = {9KG \over 3K+G}, \hspace{0.4cm}
\nu = {3K - 2G \over 2(3K+G)}.
\label{eqn:young}
\end{equation}
The effective mass density $\rho$ is given by the rule of mixtures as: $\rho = \rho_c V_c + \rho_m V_m$. The effective heat conductivity $\kappa_{\rm eff}$ and the effective coefficient of thermal expansion $\alpha_{\rm eff}$ is given by:
\begin{equation}
\frac{\kappa_{\rm eff} - \kappa_m}{\kappa_{cm}} = \frac{V_c}{1 + V_m \frac{(\kappa_c - \kappa_m)}{3\kappa_m}} \hspace{0.4cm}
\frac{\alpha_{\rm eff} - \alpha_m}{\alpha_{cm}} = \frac{ \left( \frac{1}{K_{\rm eff}} - \frac{1}{K_m} \right)}{\left(\frac{1}{K_c} - \frac{1}{K_m} \right)}
\label{eqn:thermalcondalpha}
\end{equation}
where $\kappa_{cm} = \kappa_{c} - \kappa_{m}$ and $\alpha_{cm} = \alpha_c - \alpha_m$.

\textbf{Temperature distribution through the thickness} The material properties $P$ that are temperature dependent can be written as~\cite{Reddy1998}:
\begin{equation}
P = P_o(P_{-1}T^{-1} + 1 + P_1 T + P_2 T^2 + P_3 T^3),
\end{equation}
where $P_o,P_{-1},P_1,P_2,P_3$ are the coefficients of temperature $T$ and are unique to each constituent material phase. The temperature variation is assumed to occur only in the thickness direction and the temperature field is considered to be constant in the $xy$-plane. In such a case, the temperature distribution along the thickness can be obtained by solving a steady state heat transfer problem:
\begin{equation}
-{d \over dz} \left[ \kappa(z) {dT \over dz} \right] = 0, \hspace{0.5cm} T = T_c ~\textup{at}~ z = h/2;~~ T = T_m ~\textup{at} ~z = -h/2
\label{eqn:heat}
\end{equation}
The solution of \Eref{eqn:heat} is obtained by means of a polynomial series~\cite{wu2004} as:
\begin{equation}
T(z) = T_m + (T_c - T_m) \eta(z,h)
\label{eqn:tempsolu}
\end{equation}
where,
\begin{equation}
\begin{split}
\eta(z,h) = {1 \over C} \left[ \left( {2z + h \over 2h} \right) - {\kappa_{cm} \over (k+1)\kappa_m} \left({2z + h \over 2h} \right)^{k+1} + \right. \\ 
\left. {\kappa_{cm} ^2 \over (2n+1)\kappa_m ^2 } \left({2z + h \over 2h} \right)^{2k+1}
-{\kappa_{cm} ^3 \over (3k+1)\kappa_m ^3 } \left({2z + h \over 2h} \right)^{3k+1} \right. \\ + 
\left. {\kappa_{cm} ^4 \over (4k+1)\kappa_m^4 } \left({2z + h \over 2h} \right)^{4k+1} 
- {\kappa_{cm} ^5 \over (5k+1)\kappa_m ^5 } \left({2z + h \over 2h} \right)^{5k+1} \right] ;
\end{split}
\label{eqn:heatconducres}
\end{equation}

\begin{equation}
\begin{split}
C = 1 - {\kappa_{cm} \over (k+1)\kappa_m} + {\kappa_{cm} ^2 \over (2k+1)\kappa_m ^2} 
- {\kappa_{cm} ^3 \over (3k+1)\kappa_m ^3} \\ + {\kappa_{cm} ^4 \over (4k+1)\kappa_m ^4}
- {\kappa_{cm} ^5\over (5k+1)\kappa_m ^5}
\end{split}
\end{equation}

\section{Element Formulation}
\label{fedescription}
\textbf{Shear flexible Q4 element} The plate element employed here is a $\mathcal{C}^0$ continuous shear flexible field consistent element with five degrees of freedom $(u_o,v_o,w_o,\theta_x,\theta_y)$ at four nodes in a 4-noded quadrilateral (QUAD-4) element. If the interpolation functions for QUAD-4 are used directly to interpolate the five variables $(u_o,v_o,w_o,\theta_x,\theta_y)$ in deriving the shear strains and the membrane strains, the element will lock and show oscillations in the shear and the membrane stresses. The field consistency requires that the transverse shear strains and membrane strains must be interpolated in a consistent manner. Thus, the $\theta_x$ and $\theta_y$ terms in the expressions for shear strain $\bveps_s$ have to be consistent with the derivative of the field functions, $w_{o,x}$ and $w_{o,y}$. This is achieved by using field redistributed substitute shape functions to interpolate those specific terms, which must be consistent as described in~\cite{Somashekar1987,Ganapathi1991}. This element is free from locking and has good convergence properties. For complete description of the element, interested readers are referred to the literature~\cite{Somashekar1987,Ganapathi1991}, where the element behaviour is discussed in great detail. Since the element is based on the field consistency approach, exact integration is applied for calculating various strain energy terms.

\textbf{Enriched Q4 element} In this study, an XFEM framework is employed to represent the internal discontinuity. XFEM is classified as one of the partition of unity methods. Using this property, any function $\varphi$ can be reproduced by a product of the partition of unity shape functions with $\varphi$. Let $\uu^h \subset \mathcal{U}$, the XFEM approximation can be decomposed into the standard part $\uu^h_{\rm{fem}}$ and into an enriched part $\uu^h_{\rm{enr}}$ as:

\begin{eqnarray}
\uu^h(\xx) &=& \uu^h_{\rm{fem}}(\xx) + \uu^h_{\rm{enr}}(\xx) \nonumber \\
&=& \sum\limits_{I \in \mathcal{N}^{\rm{fem}}} N_I(\xx) \qq_I + \sum\limits_{J \in \mathcal{N}^{\rm{enr}}} N_J(\xx) \varrho(\xx) \tilde{\aaa}_J
\label{eqn:xfemgeneralform}
\end{eqnarray}
where $\mathcal{N}^{\rm{fem}}$ is the set of all the nodes in the FE mesh, $\mathcal{N}^{\rm{enr}}$ is the set of nodes that are enriched with the enrichment function $\varrho$. The enrichment function $\varrho$ carries with it the nature of the solution or the information about the underlying physics of the problem, for example, $\varrho = H$, is used to capture strong discontinuities, where $H$ is the Heaviside function or $\varrho$ could be a set of functions that are chosen to represent the near tip asymptotic fields. The necessary steps involved in the implementation of the XFEM are:

\begin{enumerate}
\item {\bf Representation of the interface}  In this study, a level set approach is followed to model the cutouts. The geometric interface (for example, the boundary of the cutout) is represented by the zero level curve $\phi \equiv \phi(\xx,t) = 0$. The interface is located from the value of the level set information stored at the nodes. The standard FE shape functions can be used to interpolate $\phi$ at any point $\xx$ in the domain as:

\begin{equation}
\phi(\xx) = \sum\limits_I N_I(\xx) \phi_I
\end{equation}
where the summation is over all the nodes in the connectivity of the elements that contact $\xx$ and $\phi_I$ are the nodal values of the level set function. For circular cutout, the level set function is given by:

\begin{equation}
\phi_I = || \xx_I - \xx_c|| - r_c
\end{equation}
where $\xx_c$ and $r_c$ are the center and the radius of the cutout. For an elliptical cutout oriented at an angle $\theta$, measured from the $x-$ axis the level set function is given by:

\begin{equation}
\phi_I = \sqrt{a_1(x_I-x_c)^2-a_2(x_I-x_c)(y_I-y_c) + a_3(y_I-y_c)^2} - 1 
\end{equation}
where

\begin{equation}
a_1 = \left( {\cos \theta \over d} \right)^2, \hspace{0.15cm} a_2 = 2\cos\theta \sin\theta\left( {1 \over d^2} - {1 \over e^2} \right), \hspace{0.15cm} a_3 = \left( {\sin \theta \over d} \right)^2 + \left( {\cos \theta \over e}\right)^2.
\end{equation}
where $d$ and $e$ are the major and minor axes of the ellipse and $(x_c,y_c)$ is the center of the ellipse. \fref{fig:lvlsetinfo} shows the level set function for a circular and an elliptical cutout.

\begin{figure}[htpb]
\centering
\subfigure[Circular cutout]{\includegraphics[scale=0.6]{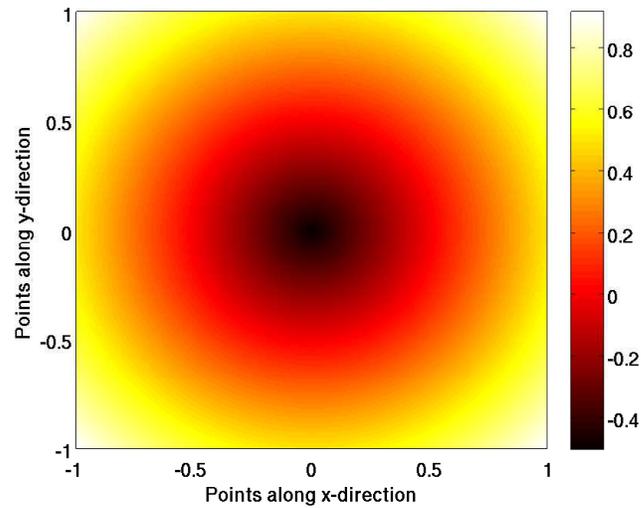}}
\subfigure[Elliptical cutout]{\includegraphics[scale=0.55]{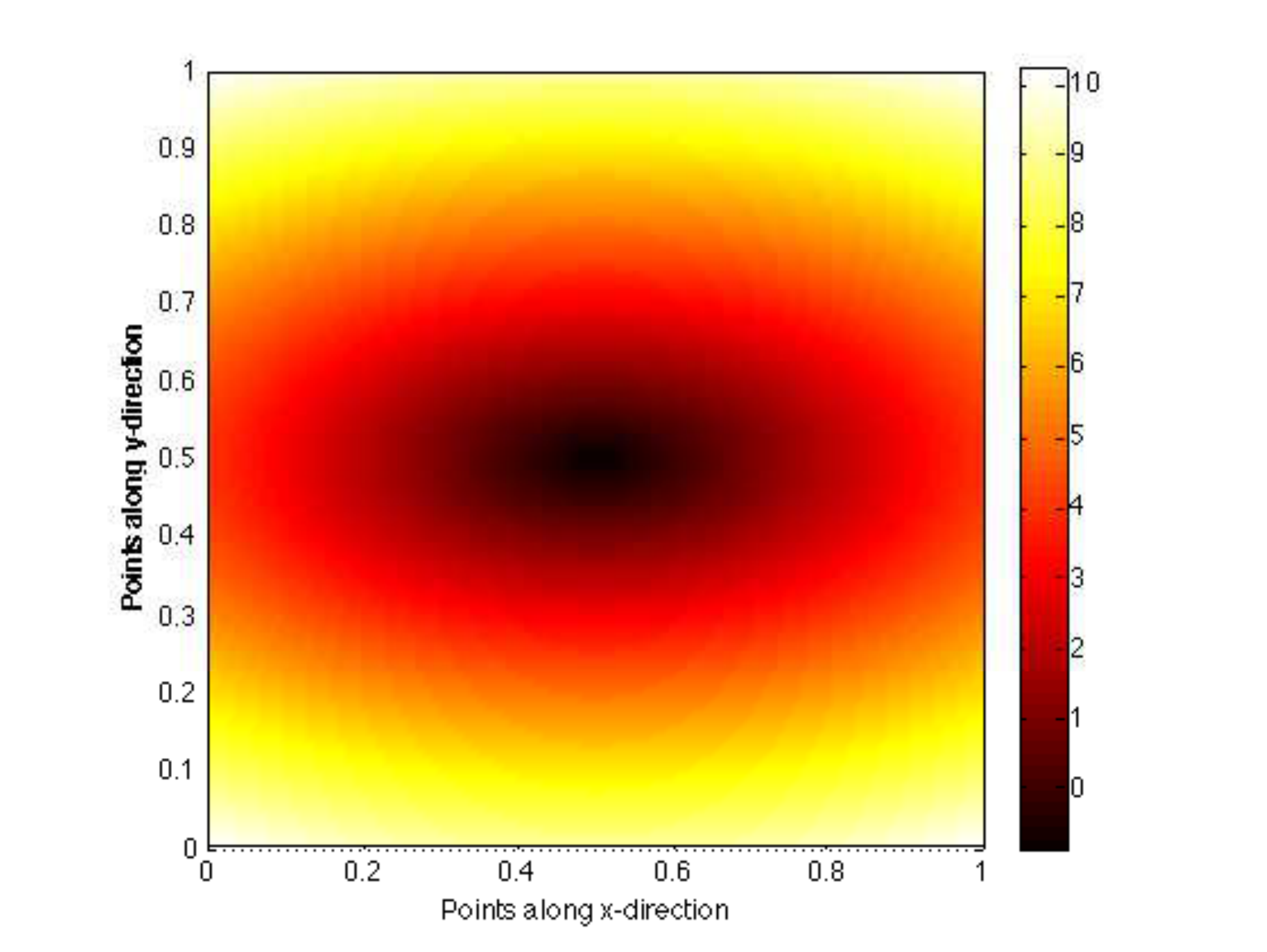}}
\caption{Level set information for different shapes.}
\label{fig:lvlsetinfo}
\end{figure}

\item {\bf Choice of enrichment function} In this study, the enriched approximation takes the following form:
\begin{equation}
\begin{split}
\uu^h(\xx)=\underbrace{ \sum\limits_{I \in \mathcal{N}^{\rm fem}} N_I(\xx) \qq_I}_{\uu^h{\rm_{std}}(\xx)}  + \\ 
\underbrace{\sum\limits_{J \in \mathcal{N}^{c}}N_J(\xx)H(\xx)\aaa_J + \sum\limits_{K \in \mathcal{N}^{f}}N_K (\xx)\sum_{\alpha=1}^4
\{\Xi_\alpha\}\bb^{\alpha}_K}_{\textup{for cracks}} +
 \underbrace{\sum\limits_{L \in \mathcal{N}^{\rm cut}} N_L (\xx) \Psi(\xx) \cc_L}_{\textup{for cutouts}}
\end{split}
\label{eqn:incrk}
\end{equation}
where $\mathcal{N}^{\rm fem}$ is the set of all nodes in the finite element mesh, $\mathcal{N}^c, \mathcal{N}^f$ and $\mathcal{N}^{\rm cut}$ are the set of nodes enriched with a Heaviside function, the near-tip asymptotic fields and the absolute value function. $\aaa_J$ and $\bb^{\alpha}_K$ are the nodal degrees of freedom corresponding to the Heaviside function $H$ and the near-tip functions, $\{\Xi_\alpha\}_{1 \leq \alpha \leq 4}$, $\cc_L$ is the nodal degrees of freedom that corresponds to the enrichment function $\Psi(\xx)$ used to represent the jump in the displacement to represent the cutouts. In this study, a Heaviside function is used to capture the physics of the problem (viz., jump across the discontinuity surface) and a set of asymptotic functions that represent the near tip fields are used, given  by:

\begin{equation}
H(\xx) = \left\{ \begin{array}{cc} 1 & \xx > 0 \\
0 & \xx < 0 \end{array} \right.
\label{eqn:heaviFn}
\end{equation}

\begin{equation}
\{\Xi_\alpha\}_{1 \le \alpha \le 4} (r,\theta)= \sqrt{r} \left\{ \sin \frac{\theta}{2}, \cos\frac{\theta}{2}, \sin\theta \sin\frac{\theta}{2},\sin\theta \cos\frac{\theta}{2} \right\}
\label{eqn:asymptotic}
\end{equation}

\item {\bf Integration}  The standard Gau\ss~quadrature cannot be applied in elements enriched by discontinuous terms, because Gau\ss~quadrature implicitly assumes a polynomial approximation. In the present study, a triangular quadrature with sub-division aligned to the discontinuous surface is employed. For the elements that are not enriched, a standard 2 $\times$ 2 Gaussian quadrature rule is used. The other techniques that can be employed are Schwarz Christoffel Mapping~\cite{natarajanbordas2009,natarajanmahapatra2010}, Generalized quadrature~\cite{mousavisukumar2011} and Smoothed eXtended FEM~\cite{bordasnatarajan2011}. 
\end{enumerate}

\begin{remark}
Appropriates terms from the displacement approximation, given by \Eref{eqn:incrk} are chosen for the study. For example, the first and the last terms are retained for the problems with cutouts only, whilst, all the terms in \Eref{eqn:incrk} are used for the problem with cracks and cutouts.
\end{remark}

\begin{remark}
The degrees of freedom of the unused nodes in the FE mesh are taken care of during the solution process.
\end{remark}

\section{Numerical Results}
\label{numerics}
In this section, we study the influence of a centrally located cutout on the fundamental frequencies of FGM plates. We consider both square and rectangular plates with simply supported and clamped boundary conditions. Two different cutout shapes, viz., circular and elliptical cutouts (see \fref{fig:problemdescription}) are considered in this study. Table \ref{table:subcellgausspt} presents the integration rules used for the current study. In all cases, we present the non dimensionalized free flexural frequencies as, unless specified otherwise:

\begin{figure}[htpb]
\centering
\includegraphics[scale=0.6]{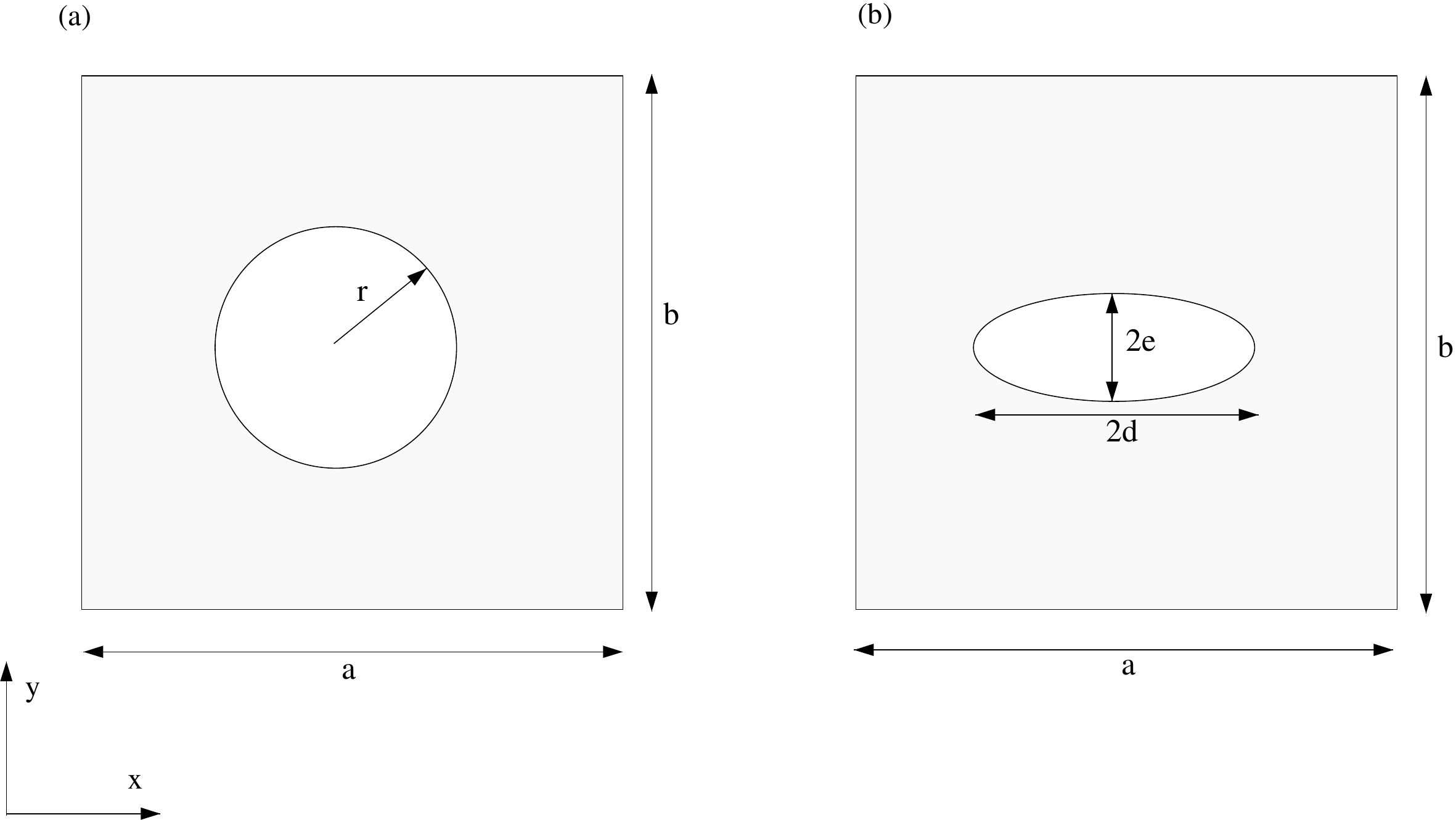}
\caption{Plate with a centrally located circular and an elliptical cutout. $r$ is the radius of the circular cutout, $2d$ and $2e$ are the major and minor axes defining the ellipse.}
\label{fig:problemdescription}
\end{figure}


\begin{table}[htpb]
\caption{Integration rules for enriched and non-enriched elements. }
\centering
\begin{tabular}{lr}
\hline
Element Type & Gau\ss ~points \\
\hline
Non-enriched element & 4  \\
Tip element & 13 per triangle  \\
Tip blending element & 16  \\
Split element & 3 per triangle  \\
Split blending element & 4  \\
Split-Tip blending element & 4 per triangle \\
\hline
\end{tabular}
\label{table:subcellgausspt}
\end{table}

\begin{equation}
\Omega = \omega \left( \frac{a^2}{h} \right) \sqrt{ \frac{\rho_c h}{D_c}}
\label{eqn:nondimfreq}
\end{equation}
where $D_c = E_ch^3/(12(1-\nu^2))$ and $\rho_c$ is the mass density. In order to be consistent with the existing literature, properties of the ceramic are used for normalization. The effect of the plate slenderness ratio $a/h$, the  plate aspect ratio $b/a$, the cutout radius $r/a$, the cutout geometry $d/e$ and the boundary condition on the natural frequencies are numerically studied. The FGM plate considered here consists of silicon nitride (Si$_3$N$_4$) and stainless steel (SUS304). The material is considered to be temperature dependent and the temperature coefficients corresponding to Si$_3$N$_4$/SUS304 are listed in Table \ref{table:tempdepprop} ~\cite{Sundararajan2005,Reddy1998}. The mass density $(\rho)$ and the thermal conductivity $(K)$ are: $\rho_c=$ 2370 kg/m$^3$, $K_c=$ 9.19 W/mK for Si$_3$N$_4$ and $\rho_m =$ 8166 kg/m$^3$, $K_m =$ 12.04 W/mK for SUS304. Poisson's ratio $\nu$ is assumed to be constant and taken as 0.28 for the current study~\cite{Sundararajan2005,Prakash2007}.  Here, the modified shear correction factor obtained based on energy equivalence principle as outlined in~\cite{Singh2011} is used. The boundary conditions for simply supported and clamped cases are :

\noindent \emph{Simply supported boundary condition}: \\
\begin{equation}
u_o = w_o = \theta_y = 0 \hspace{0.2cm} ~\textup{on} \hspace{0.2cm}  x=0,a; \hspace{0.2cm}
v_o = w_o = \theta_x = 0 \hspace{0.2cm} ~\textup{on} \hspace{0.2cm}  y=0,b
\end{equation}

\noindent \emph{Clamped boundary condition}: \\
\begin{equation}
u_o = w_o = \theta_y = v_o = \theta_x = 0 \hspace{1cm} ~\textup{on} ~ x=0,a \hspace{0.2cm} \& \hspace{0.2cm}  y=0,b
\end{equation}

\begin{table}
\renewcommand\arraystretch{1.5}
\caption{Temperature dependent coefficient for material Si$_3$N$_4$/SUS304, Ref~\cite{Reddy1998,Sundararajan2005}.}
\centering
\begin{tabular}{lcccccc}
\hline
Material & Property & $P_o$ & $P_{-1}$ & $P_1$ & $P_2$ & $P_3$  \\
\hline
\multirow{2}{*}{Si$_3$N$_4$} & $E$(Pa) & 348.43e$^9$ &0.0& -3.070e$^{-4}$ & 2.160e$^{-7}$ & -8.946$e^{-11}$  \\
& $\alpha$ (1/K) & 5.8723e$^{-6}$ & 0.0 & 9.095e$^{-4}$ & 0.0 & 0.0 \\
\cline{2-7}
\multirow{2}{*}{SUS304} & $E$(Pa) & 201.04e$^9$ &0.0& 3.079e$^{-4}$ & -6.534e$^{-7}$ & 0.0  \\
& $\alpha$ (1/K) & 12.330e$^{-6}$ & 0.0 & 8.086e$^{-4}$ & 0.0 & 0.0 \\
\hline
\end{tabular}
\label{table:tempdepprop}
\end{table}

\paragraph{Validation} Before proceeding with a detailed study on the effect of different parameters on the natural frequency, the formulation developed herein is validated against available results pertaining to the linear frequencies of a FGM plate in thermal environment~\cite{Huang2004} and an isotropic plate with a centrally located circular and an elliptical cutout~\cite{huangsakiyama1999}. The computed linear frequencies: (a) for a square simply supported FGM plate in thermal environment with $a/h=$ 8 is given in Table \ref{table:tempValidation} and (b) for square plate with a circular cutout for various boundary conditions is given in Table \ref{table:cutoutvalida}. It can be seen that the numerical results from the present formulation are found to be in good agreement with the existing solutions. Based on a progressive refinement a 40 $\times$ 40 quadrilateral mesh is found to be adequate to model the full plate with a cutout.
\begin{table}
\renewcommand\arraystretch{1}
\caption{Comparison of non-dimensional linear frequency $\left( \Omega = \omega \left( \frac{a^2}{h} \right) \sqrt{ \frac{\rho_m (1-\nu^2)}{E_m}} \right)$ of simply supported FGM plate $(a/b=1, a/h=8)$ in thermal environment.}
\centering
\begin{tabular}{ccrrrrr}
\hline
Temperature & gradient index & \multicolumn{2}{c}{Mode 1} & & \multicolumn{2}{c}{Mode 2} \\
\cline{3-4} \cline{5-7}
$T_c,T_m$& $k$ & Ref.~\cite{Huang2004} & Present & & Ref.~\cite{Huang2004} & Present \\
\hline
\multirow{4}{*}{$T_c=$ 400K, $T_m=$ 300K} & 0.0 & 12.397 & 12.315 & & 29.083 & 29.031 \\
& 0.5 & 8.615 & 8.484 && 20.215 & 19.986 \\
& 1.0 & 7.474 & 7.443 & & 17.607 & 17.515 \\
& 2.0 & 6.693 & 6.679 & & 15.762 & 15.709 \\
\cline{2-7}
\multirow{4 }{*}{$T_c=$ 600K, $T_m=$ 300K} & 0.0 & 11.984 & 11.894 & & 28.504 & 28.436 \\
& 0.5 & 8.269 & 8.147 & & 19.784 & 19.535 \\
& 1.0 & 7.171 & 7.126 & & 17.213 & 17.098 \\
& 2.0 & 6.398 & 6.370 & & 15.384 & 15.309 \\
\hline
\end{tabular}
\label{table:tempValidation}
\end{table}

\begin{table}
\renewcommand\arraystretch{1}
\caption{Natural frequency parameter $\left( \Omega = \left[ \frac{ \omega^2 \rho_c h a^4}{D_c (1-\nu^2)} \right]^{1/4} \right)$ for an isotropic plate with a central cutout.}
\centering
\begin{tabular}{lccrr}
\hline
Cutout & Boundary & Mesh & \multicolumn{2}{c}{Frequency} \\
\cline{4-5}
 &  Condition & & mode 1 & mode 2 \\
\hline
\multirow{4}{*}{Circular ($r/a=$ 0.1)} & \multirow{4}{*}{CCCC} & 20$\times$20 & 6.1848 & 8.7215 \\
& & 30$\times$30 & 6.1762 & 8.6622 \\
& & 40$\times$40 & 6.1725 & 8.6443 \\
& & Ref.~\cite{huangsakiyama1999} & 6.2110 & 8.7310 \\
\cline{2-5}
\multirow{4}{*}{Elliptical, ($2d/a=$3/8; $2d/2e=$ 2)} & \multirow{4}{*}{SSSS} & 20 $\times$20 & 4.4828 & 6.9237 \\
& & 30$\times$30 & 4.4775 & 6.8849\\
& & 40$\times$40 & 4.4758 & 6.8705 \\
& & Ref.~\cite{huangsakiyama1999} & 4.4820 & 6.9120 \\
\hline
\end{tabular}
\label{table:cutoutvalida}
\end{table}

\begin{figure}[htpb]
\centering
\includegraphics[scale=0.6]{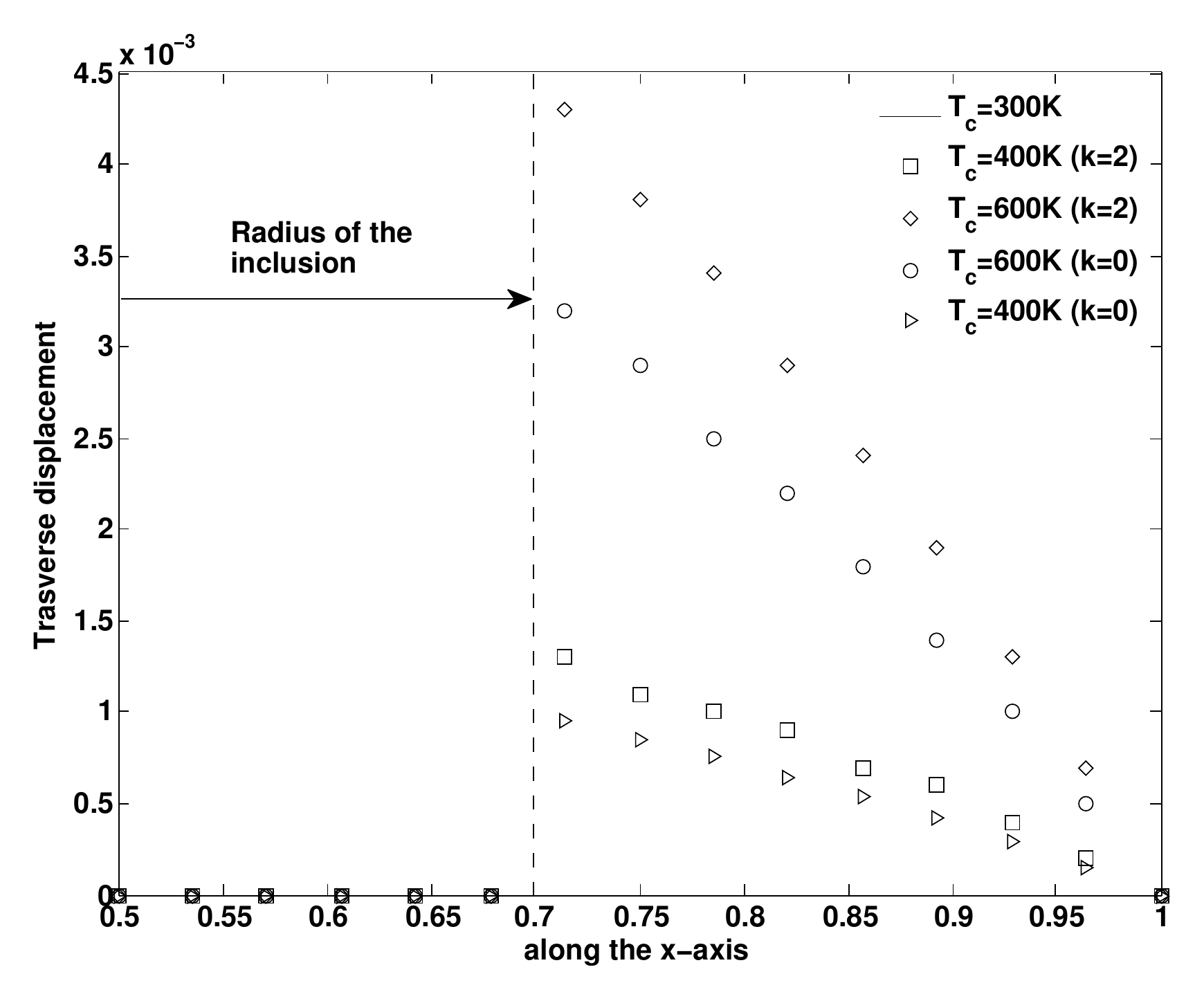}
\caption{Transverse displacement along the centerline $(w(x,b/2,0))$ of the plate for various gradient indices ($k=$ 0, 2) and for different thermal gradient, viz., $\Delta T=$ 0K, 100K, 300K ($T_c=$ 300K, 400K and 600K and $T_m=$ 300K) with a circular cutout $(r/a=0.2)$ at the center.  }
\label{fig:thermDisplace}
\end{figure}

\begin{figure}[htpb]
\centering
\includegraphics[scale=0.6]{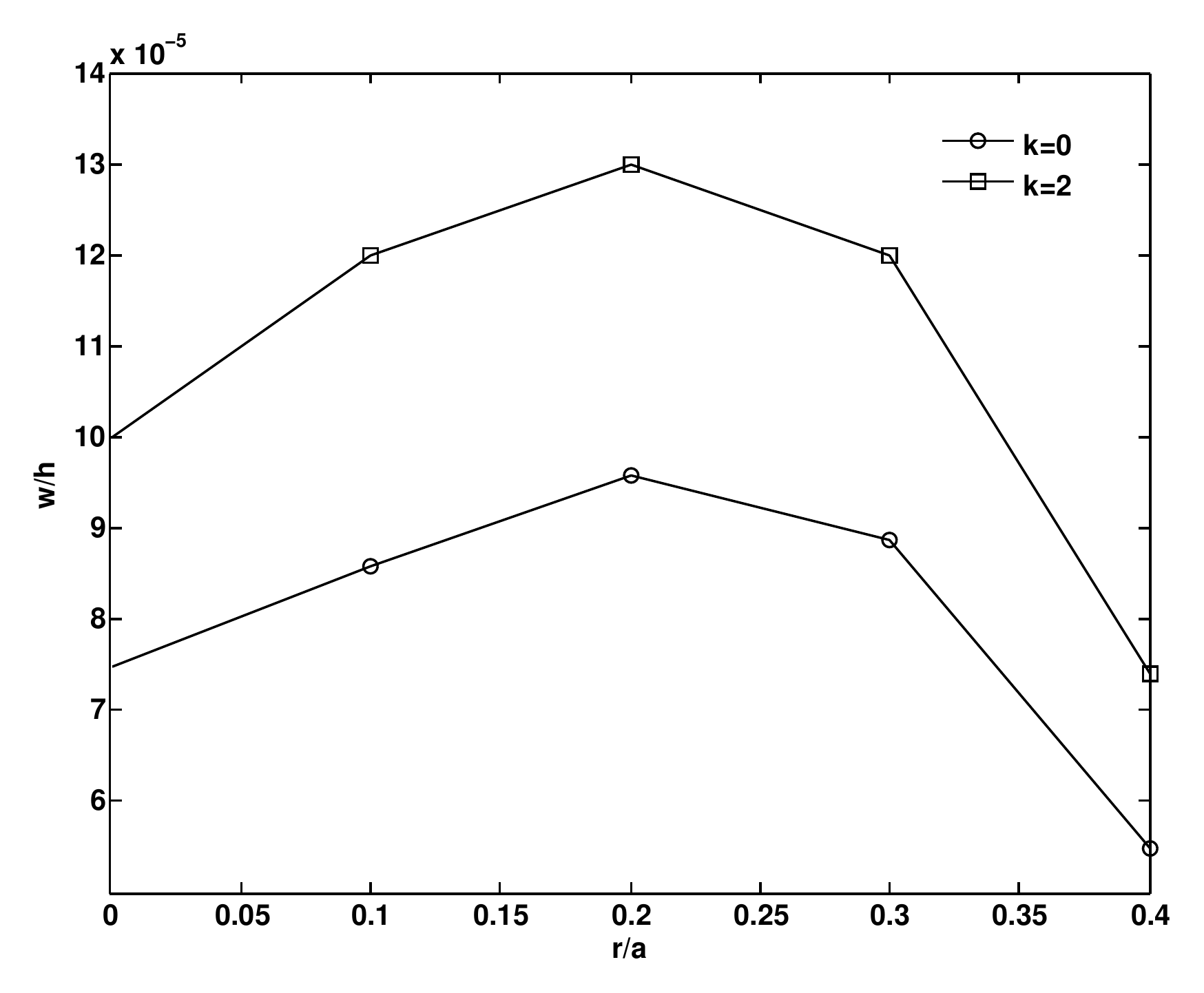}
\caption{Maximum transverse displacement $(w_{max} = w/h)$ as a function of the cutout radius for a square simply supported FGM plate for different gradient index ($k=0,2$) subjected to a thermal gradient $\Delta T = 100K$.}
\label{fig:RtoAWmax}
\end{figure}

\paragraph{Effect of thermal environment}
Next, the linear free flexural vibration behaviour of FGM is numerically studied with and without thermal environment. For the uniform temperature case, the material properties are evaluated at $T_c = T_m = 300K$. The temperature field is assumed to vary only in the thickness direction and determined by \Eref{eqn:tempsolu}. The temperature for the ceramic surface is varied $(T_c = 400K, 600K, 900K)$, whilst maintaining a constant value on the metallic surface $(T_m=300K)$ to subject the plate to a thermal gradient. \fref{fig:thermDisplace} shows the transverse displacement along $y=b/2$ for a simply supported square FGM plate under different thermal loading with center circular cutout $r/a=$ 0.2. It can be seen that the transverse displacement increases with increasing thermal gradient. The transverse displacement for a FGM plate with gradient index $k=$ 2 is greater than a plate with $k=$ 0 (note that $k=$ 0 is a pure ceramic plate) as expected. \fref{fig:RtoAWmax} shows the influence of cutout radius on the maximum transverse deflection for a simply supported square FGM plate in thermal environment, $\Delta T = 100K$. It can be seen that increasing the gradient index, increases the transverse deflection. With increase in the cutout radius, the maximum transverse displacement first increases due to stiffness reduction and with further increase the transverse displacement decreases. This is because for a larger cutout radius, there is less material to deform. The increase in the transverse displacement with increasing gradient index can be attributed to the increase in the metallic volume fraction, which increases the coefficient of thermal expansion and thus increasing the transverse displacement. The increase in the displacement with cutout radius is due to stiffness degradation.

\begin{figure}[htpb]
\centering
\includegraphics[scale=0.6,angle=90]{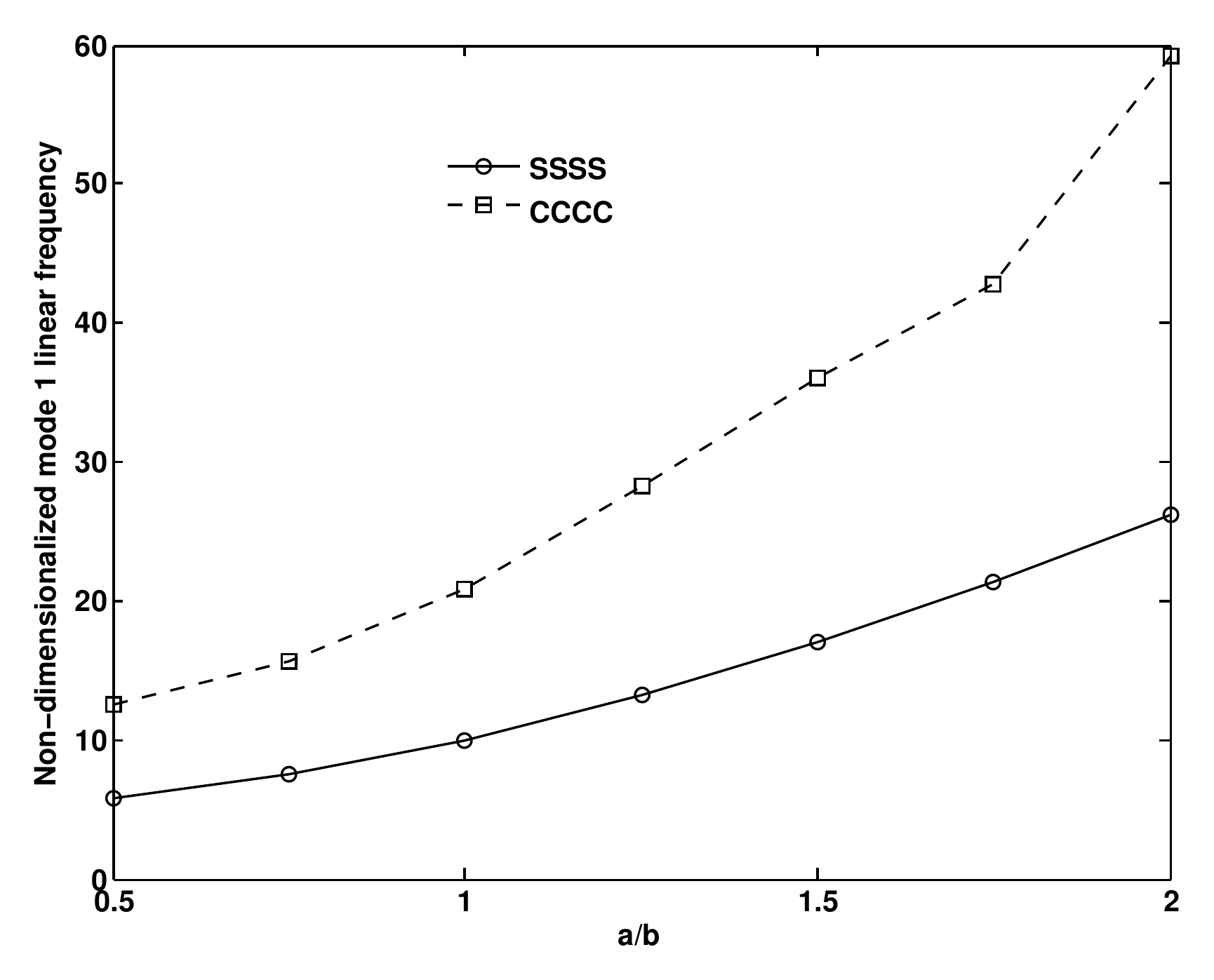}
\caption{Effect of the plate aspect ratio on the linear frequency for a FGM plate with $a/h=$10 and with a thermal gradient $\Delta T=$ 100K ($T_c=$ 400K, $T_m=$ 300K) for a centrally located circular cutout with radius $r/a=$ 0.2 for various boundary conditions (simply supported and clamped) with gradient index $k=$ 2. Note that in this study, $a$ is kept constant.}
\label{fig:effplateaspect}
\end{figure}

The geometric stiffness matrix is computed from the in-plane stress resultants due to the applied thermal gradient. The geometric stiffness matrix is then added to the stiffness matrix and then the eigenvalue problem is solved. Table \ref{table:platethicknesseffect} shows the influence of the gradient index, the thermal gradient and the plate aspect ratio on the fundamental frequency of a square simply supported FGM plate with a circular cutout $r/a=0.2$. It can be seen that the combined effect of increasing the gradient index and the thermal gradient is to lower the fundamental frequency, whilst the frequency increases with decreasing plate thickness. The effect of boundary conditions, the plate thickness and the gradient index $k$ on the fundamental frequency for a square FGM plate in thermal environment ($\Delta T=100K$) is shown in Table \ref{table:platebcffect}. The fundamental frequency decreases with increasing gradient index due to increase in the metallic volume fraction. The frequency initially increases with plate thickness but upon further increase, the fundamental frequency decreases. \fref{fig:effplateaspect} shows the influence of the plate aspect ratio on the fundamental frequency for a FGM plate in thermal environment $\Delta T = 100K$ for different boundary conditions with gradient index $k=$ 2 and $a/h=$ 10. It can be seen that the frequency increases with increasing $a/b$ ratio and clamped plate has higher frequency than a simply supported plate. The increase in the stiffness is the cause for increase in frequency when the boundary condition is changed from simply supported to clamped condition for a fixed aspect ratio and plate thickness.

\begin{table}
\renewcommand\arraystretch{1.5}
\caption{Non-dimensionalized mode 1 frequency for square plates with circular cutouts $r/a=$ 0.2 with $a/h=$ 5, 10 in thermal environment.}
\centering
\begin{tabular}{ccrrrr}
\hline
$a/h$ & $k$ & $T_c=$ 300K & $T_c=$ 400K & $T_c=$ 600K & $T_c=$ 900K\\
\hline
\multirow{5}{*}{5} & 0 & 17.6855 & 17.4690 & 17.0266 & 16.3111\\
& 1 & 10.6681 & 10.5174 & 10.1932 & 9.6350 \\
& 2 & 9.6040 & 9.4618 & 9.1469 & 8.5882\\
& 5 & 8.7113 & 8.5738 & 8.2544 & 7.6601 \\
& 10 & 8.2850 & 8.1484 & 7.8191 & 7.1840 \\
\cline{2-6}
\multirow{5}{*}{10} & 0 & 19.1844 & 18.5992 & 17.2928 & 14.8701 \\
& 1 & 11.5736 & 11.1317 & 10.1161 & 8.1589 \\
& 2 & 10.4135 & 9.9844 & 8.9842 & 7.0277\\
& 5 & 9.4461 & 9.0145 & 7.9853 & 5.9172 \\
& 10 & 8.9858 & 8.5452 & 7.4747 & 5.2682 \\
\hline
\end{tabular}
\label{table:platethicknesseffect}
\end{table}

\begin{table}
\renewcommand\arraystretch{1.5}
\caption{Non-dimensionalized mode 1 linear frequency for square plates with circular cutouts $r/a=$ 0.2, temperature gradient $\Delta T=$ 100 K ($T_c=$ 400K, $T_m=$ 300K).}
\centering
\begin{tabular}{ccrrrrr}
\hline
Boundary & $a/h$ & \multicolumn{3}{c}{gradient index, $k$} \\
\cline{3-7}
Condition & & 0 & 1 & 2 & 5 & 10\\
\hline
\multirow{4}{*}{SSSS} & 5 &17.4690 & 10.5174 & 9.4618 & 8.5738 & 8.1484 \\
& 10 & 18.5992 & 11.1317 & 9.9844 & 9.0145 & 8.5452\\
& 20 & 17.5380 & 10.1776 & 9.0071 & 7.9664 & 7.4294 \\
& 25 & 16.3587 & 9.1762 & 7.9919 & 6.8810 & 6.2689\\
\cline{2-7}
\multirow{4}{*}{CCCC} & 5 & 31.4944 & 18.9259 & 16.9575 & 15.3461 & 14.6221 \\
& 10 & 38.7777 & 23.2411 & 20.8571 & 18.9072 & 18.0010\\
& 20 & 41.0541 & 24.4016 & 21.8525 & 19.7479 & 18.7389 \\
& 25 & 40.7846 & 24.0796 & 21.5110 & 19.3700 & 18.3275\\
\hline
\end{tabular}
\label{table:platebcffect}
\end{table}

\begin{figure}[htpb]
\centering
\includegraphics[scale=0.6]{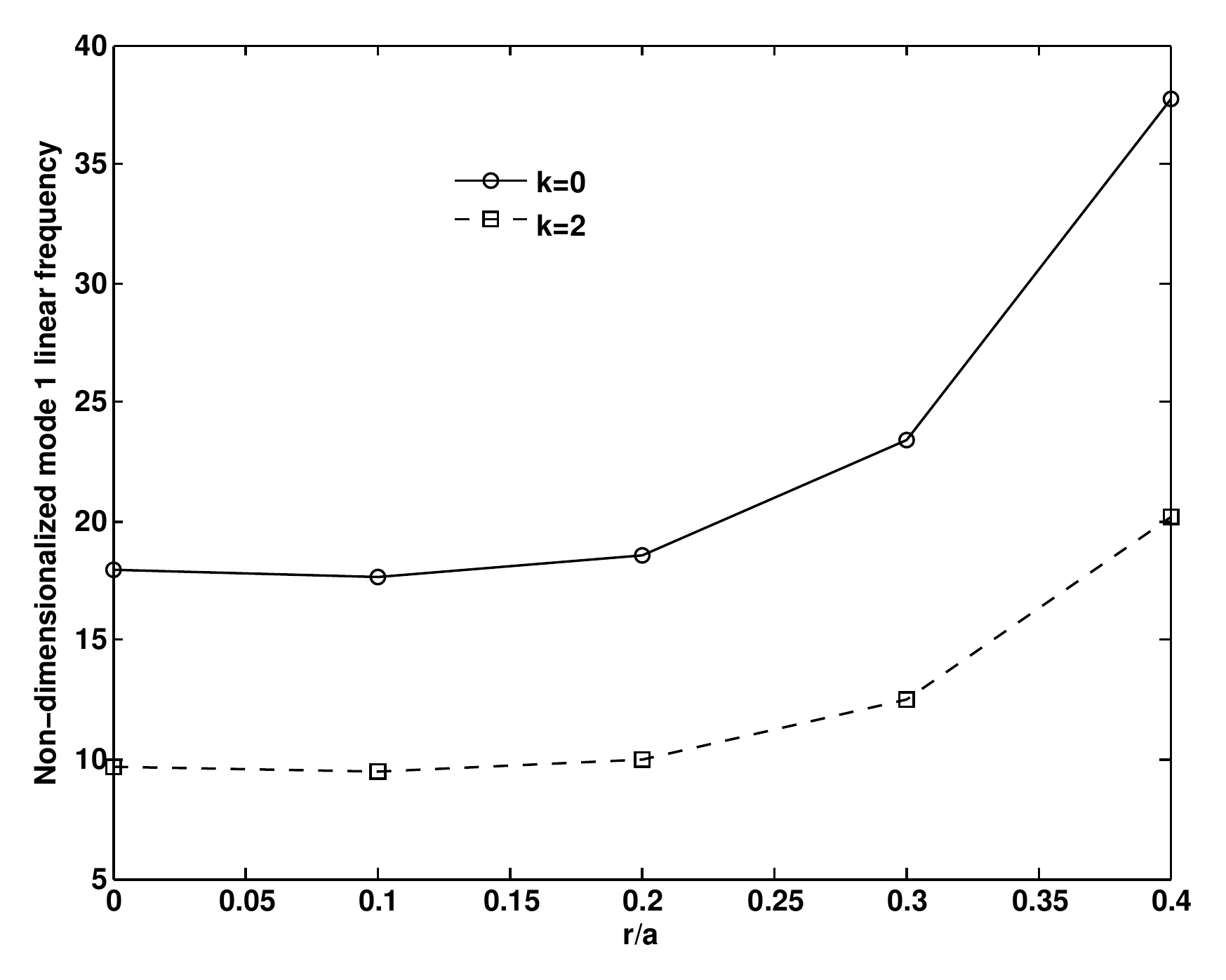}
\caption{Effect of the cutout size on the linear frequency $(\Omega)$ for a square simply supported FGM plate with $a/h=$10 in a thermal environment $\Delta T=$ 100K ($T_c=$ 400K, $T_m=$ 300K) for a centrally located cutout for various gradient indices $k$.}
\label{fig:effcutoutsize}
\end{figure}

\begin{figure}[htpb]
\centering
\includegraphics[scale=0.6]{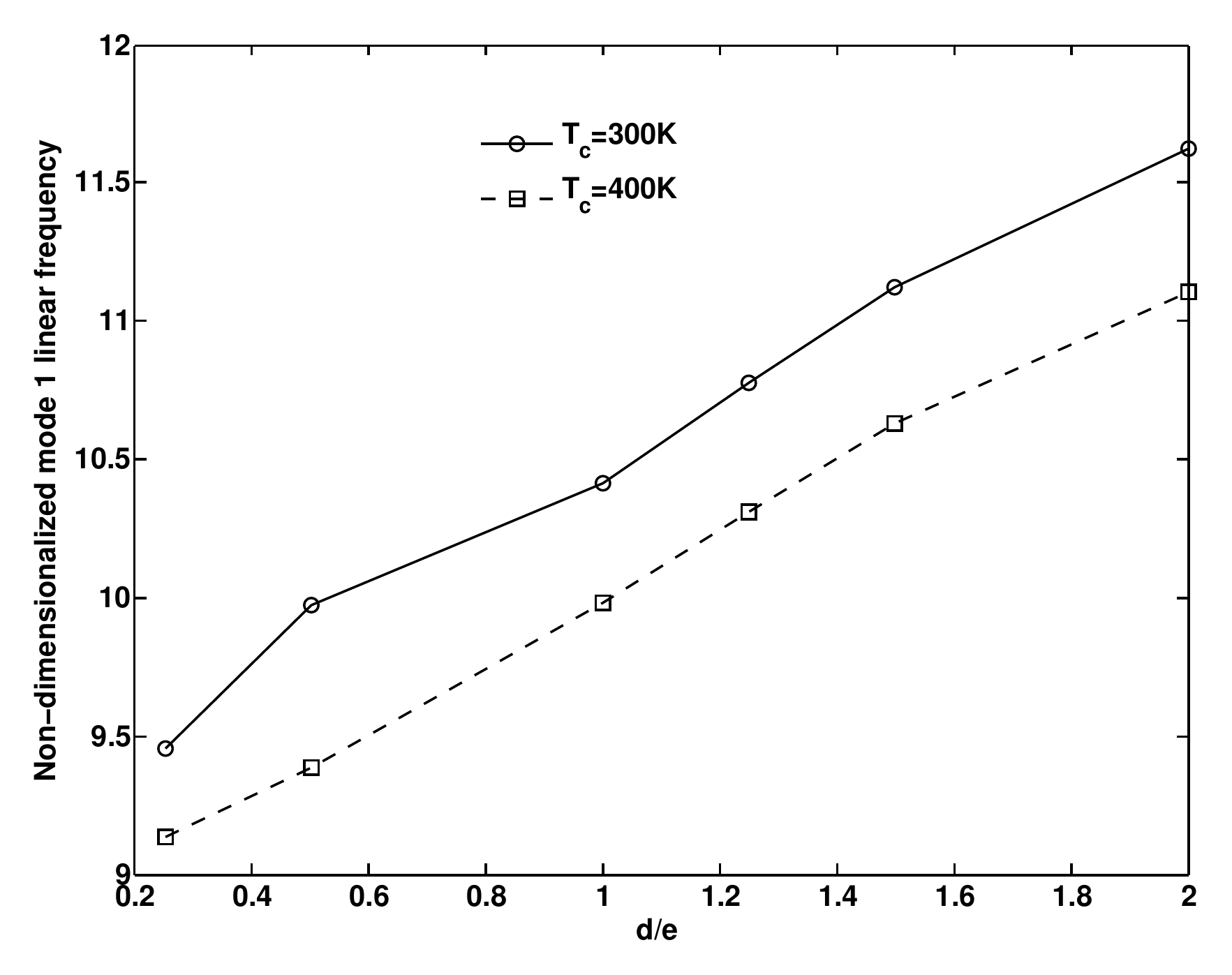}
\caption{Influence of the geometry of the cutout on the linear frequency $(\Omega)$ for a square FGM plate with $a/h=$10 with gradient index $k=$ 2 for various thermal gradients.}
\label{fig:effcutoutgeom}
\end{figure}

\paragraph{Effect of cutout geometry}
\fref{fig:effcutoutsize} shows the influence of the cutout size on the frequency for a plate in thermal environment ($\Delta T=100K$). The frequency increases with increasing cutout size, whilst decreases with increasing gradient index. The effect of geometry of the cutout $d/e$ is shown in \fref{fig:effcutoutgeom} for a square FGM plate in two different thermal environment, viz., $\Delta T = 0K, 100K$. It can be seen that the fundamental frequency increases with increasing the $d/e$ ratio, whilst decreases with increasing thermal gradient. 

\begin{figure}[htpb]
\centering
\includegraphics[scale=0.8]{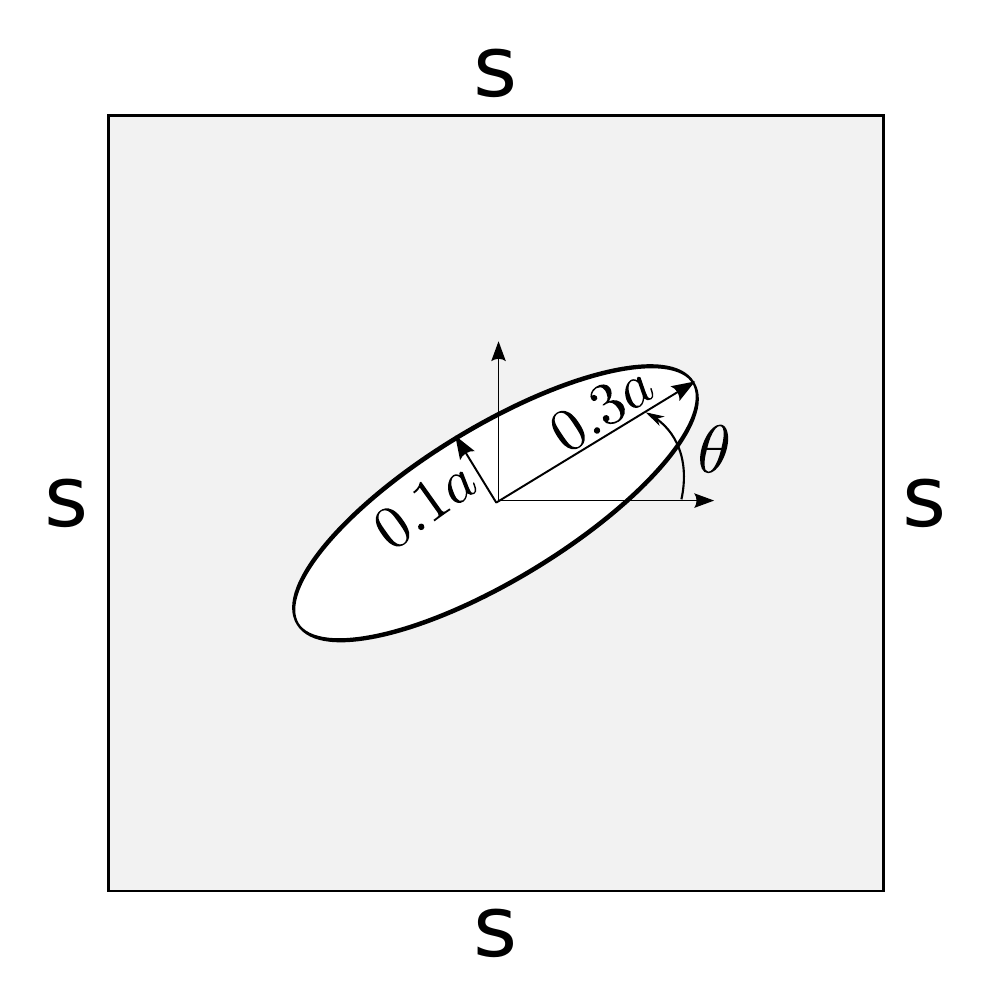}
\caption{Plate with a centrally located elliptical cutout oriented at an angle $\theta$, measured in the counterwise direction from the $x-$ axis.}
\label{fig:inclinedEllipse}
\end{figure}

Table \ref{table:Elippiticcutoutorientation} shows the influence of the orientation $\theta$ of an elliptical cutout on the fundamental frequencies for a square simply supported FGM plate in thermal environment with $a/h=$ 10 and gradient index $k=$ 2. The geometry of the plate is shown in \fref{fig:inclinedEllipse}. It can be seen that the orientation $\theta$ of the cutout has strong influence on the fundamental frequencies. It can be seen that increasing the orientation $\theta$ from 0$^\circ$ to 90$^\circ$, the mode 1 and the mode 3 first decreases until $\theta=$ 45$^\circ$ and with further increase in $\theta$, the mode 1 and mode 3 frequency increases. On the other hand, the mode 2 frequency first increases until $\theta=$ 45$^\circ$ and then decreases. The frequencies are symmetric with respect to $\theta=$ 45$^\circ$, indicating that the FGM plate is globally homogeneous. With increasing temperature, all the fundamental frequencies decreases as observed earlier.


\begin{table}
\renewcommand\arraystretch{1.5}
\caption{Non-dimensionalized natural frequency for a square simply supported FGM plate ($a/b=$ 1, $a/h=$ 10, $k=$ 2) with an elliptical cutout.}
\centering
\begin{tabular}{c r r r r r r}
\hline
Major axis &\multicolumn{3}{c}{$T_c=$ 300K}	   &\multicolumn{3}{c}{$T_c=$ 400K} \\
\cline{2-4} \cline{5-7} 

  orientation 	& mode 1	& mode 2	& mode 3 	& mode 1	& mode 2	& mode 3  \\
\hline

  0$^\circ$	&9.6795  & 16.3035  & 23.0825   & 9.2771  & 16.0082  & 22.7147 \\
 10$^\circ$	&9.6627  & 16.3067  & 23.0677   & 9.2591  & 16.0108  & 22.6994 \\
 20$^\circ$	&9.6231  & 16.3311  & 23.0329   & 9.2168  & 16.0338  & 22.6632 \\
 30$^\circ$	&9.5740  & 16.3353  & 22.9892   & 9.1643  & 16.0364  & 22.6175 \\
 40$^\circ$	&9.5435  & 16.3548  & 22.9596   & 9.1317  & 16.0548  & 22.5868 \\
 50$^\circ$	&9.5435  & 16.3548  & 22.9596   & 9.1317  & 16.0548  & 22.5868 \\
 60$^\circ$	&9.5740  & 16.3353  & 22.9892   & 9.1643  & 16.0364  & 22.6175 \\
 70$^\circ$	&9.6231  & 16.3311  & 23.0329   & 9.2168  & 16.0338  & 22.6632 \\
 80$^\circ$	&9.6627  & 16.3067  & 23.0677   & 9.2591  & 16.0108  & 22.6994 \\
 90$^\circ$	&9.6795  & 16.3035  & 23.0825   & 9.2771  & 16.0082  & 22.7147 \\

\hline
\label{table:Elippiticcutoutorientation}
\end{tabular}
\end{table}

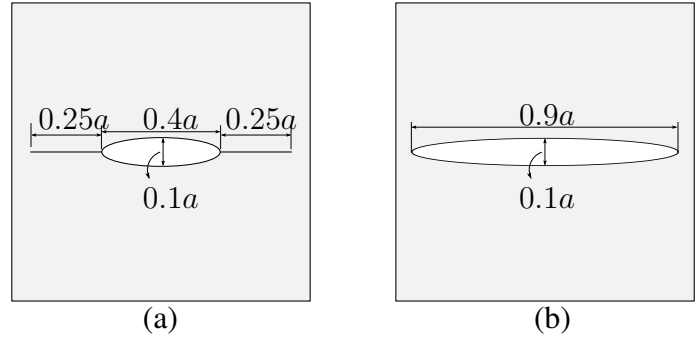
\begin{figure}[ht!]
\centering
\begin{overpic}[scale=.5,unit=1mm]%
{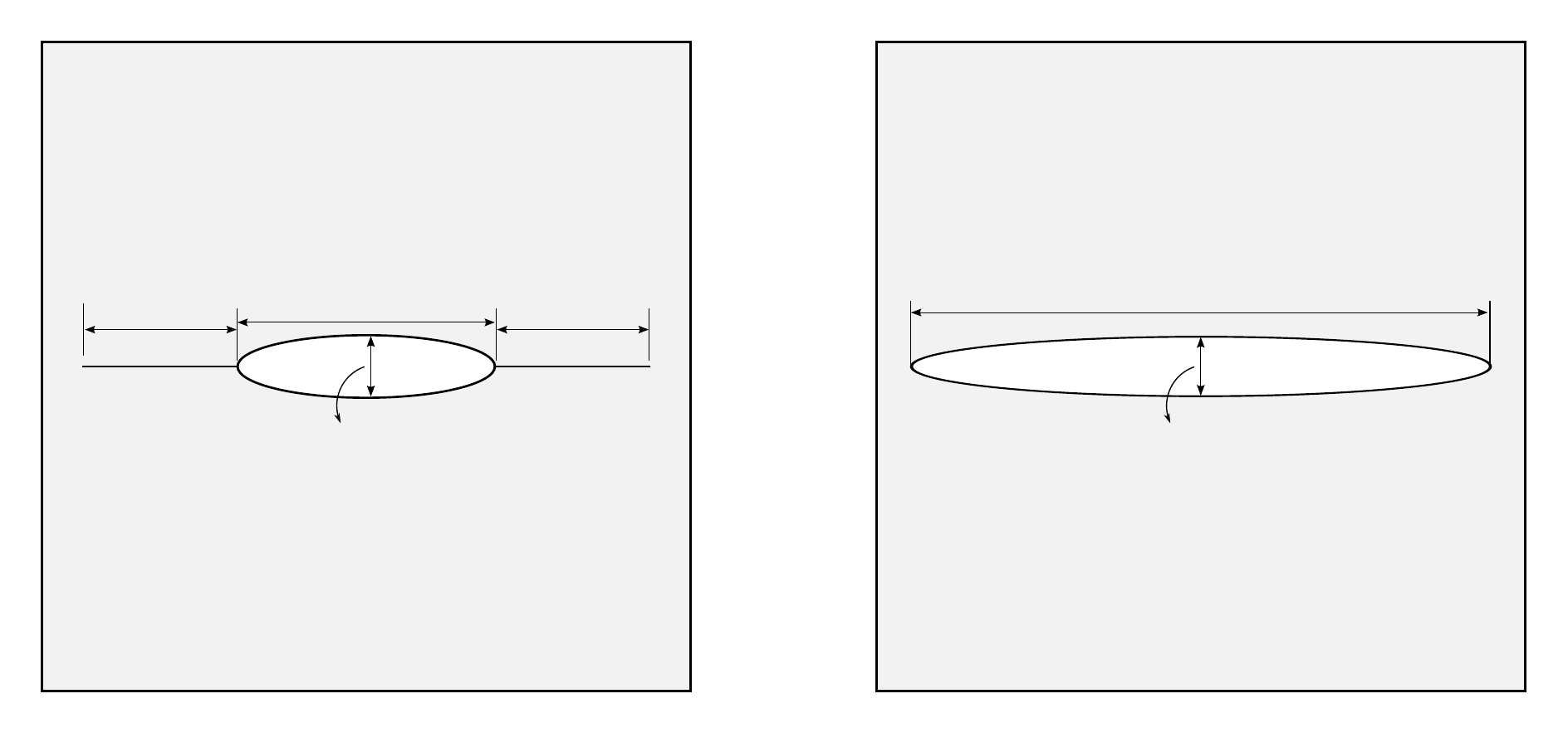}
\put(20,-1){(a)}
\put(72,-1){(b)}
\put(6.1,25.5){$0.25a$}
\put(31,25.5){$0.25a$}
\put(70,26){$0.9a$}
\put(20,25.5){$0.4a$}
\put(20,15){$0.1a$}
\put(70,15){$0.1a$}
\end{overpic}
\caption{Simply supported FGM plates $a/h=10$,  a) with an elliptic cutout and two cracks b) with an elliptic crack.}
\label{fig:EllipsCutoutCracks}
\end{figure}

\paragraph{Cracks emanating from the cutout}
Next, we study the influence of cracks emanating from an elliptical cutout. Two cases are considered as show in \fref{fig:EllipsCutoutCracks}, viz., (a) Case A: elliptical cutout with two cracks one on either side with crack length $\ell=$ 0.25$a$ (b) Case B: single elliptical cutout with length of the major axis equals the sum of the lengths of the two cracks and an elliptical cutout in Case A. Table \ref{table:cutoutwithcrack} shows the influence of the gradient index and the thermal gradient on the natural frequencies of a square simply supported FGM plates. It can be seen that increasing the temperature of the ceramic phase and the gradient index, the fundamental frequency decreases. The mode 1 and mode 3 frequency for both the cases show very similar behaviour, whilst the mode 2 frequency for Case A is greater than the mode 2 frequency for Case B. This can be attributed to the fact that Case A can considerably more material to resist the deformation and that in the case of mode 2, the crack is parallel to the deformed shape. The mode 2 deformed shape can be seen as the crack opening mode. On the other hand, the crack and the cutout is perpendicular to the deformed shape in the case of mode 3, which resists the crack opening.


\begin{table}
\renewcommand\arraystretch{1.5}
\caption{Non-dimensionalized natural frequency for square FGM plates with: (a) elliptic cutout with two cracks and (b) elliptic crack }
\centering
\begin{tabular}{c c r r r r r r}
\hline
			& 	  &\multicolumn{3}{c}{ Case A: ellipse with two cracks}	   &\multicolumn{3}{c}{Case B: elliptic crack} \\
\cline{3-5} \cline{6-8} 

$\Delta T$		& k	& mode 1	& mode 2	& mode 3 	& mode 1	& mode 2	& mode 3  \\
\cline{1-5} \cline{6-8} 

\multirow{3}*{0}
			& 0	& 16.2790	& 34.4610	& 41.5259	& 16.1708  	& 20.4605  	& 42.7517 \\
			& 2	&  8.8710	& 18.6648	& 22.5162	& 8.8104   	& 11.1051  	& 23.1853\\
			& 5	&  8.0473	& 16.8857	& 20.4443	& 7.9928  	& 10.0659  	& 21.0584\\
\cline{2-5} \cline{6-8} 

\multirow{3}*{100}
			& 0	&15.7865	& 34.0478	& 40.9831	& 15.6527	& 20.0172 	& 42.1728 \\
			& 2	& 8.5066	& 18.3994	& 22.1638	&  8.4238 	& 10.7881  	& 22.8061 \\
			& 5	& 7.6827	& 16.6294	& 20.0996	&  7.6057 	&  9.7508  	& 20.6869 \\
\hline
\label{table:cutoutwithcrack}
\end{tabular}
\end{table}

\begin{figure}[ht!]
\centering
\includegraphics[scale=.75]{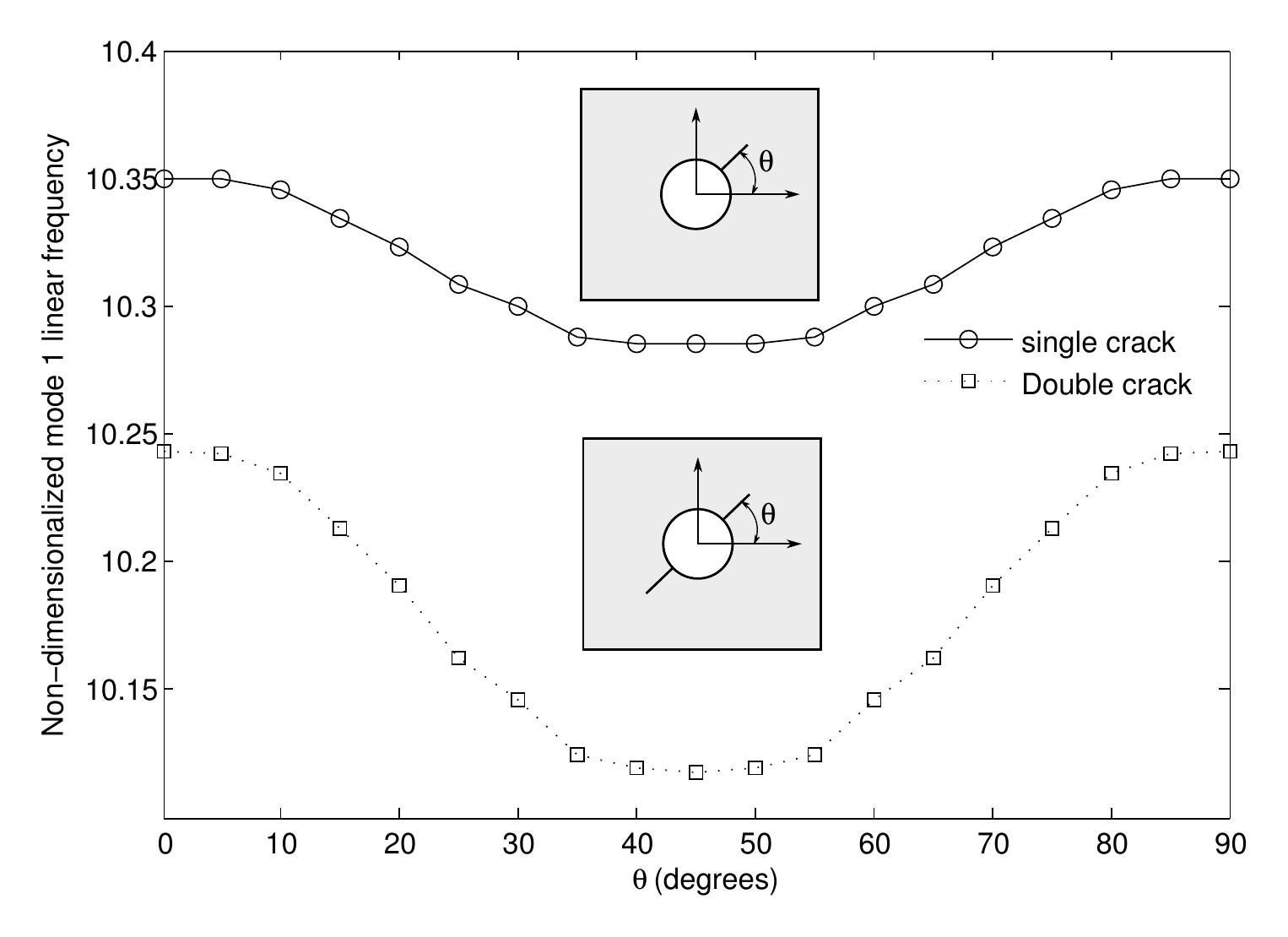}
\caption{The effect of the crack orientation $\theta$ on the mode 1 natural frequency $(\Omega)$ for a square simply supported FGM plate ($a/h=$10, $k=$2) with a circular cutout $r/a=0.2$ and $\Delta T = $0 ($T_c=T_m=$ 300K). The length of each crack is $\ell/a=0.1$.}
\label{fig:NF_c1c2_n2_T300_Rc01_ssss}
\end{figure}

\begin{figure}[ht!]
\centering
\begin{overpic}[scale=.75,unit=1mm]{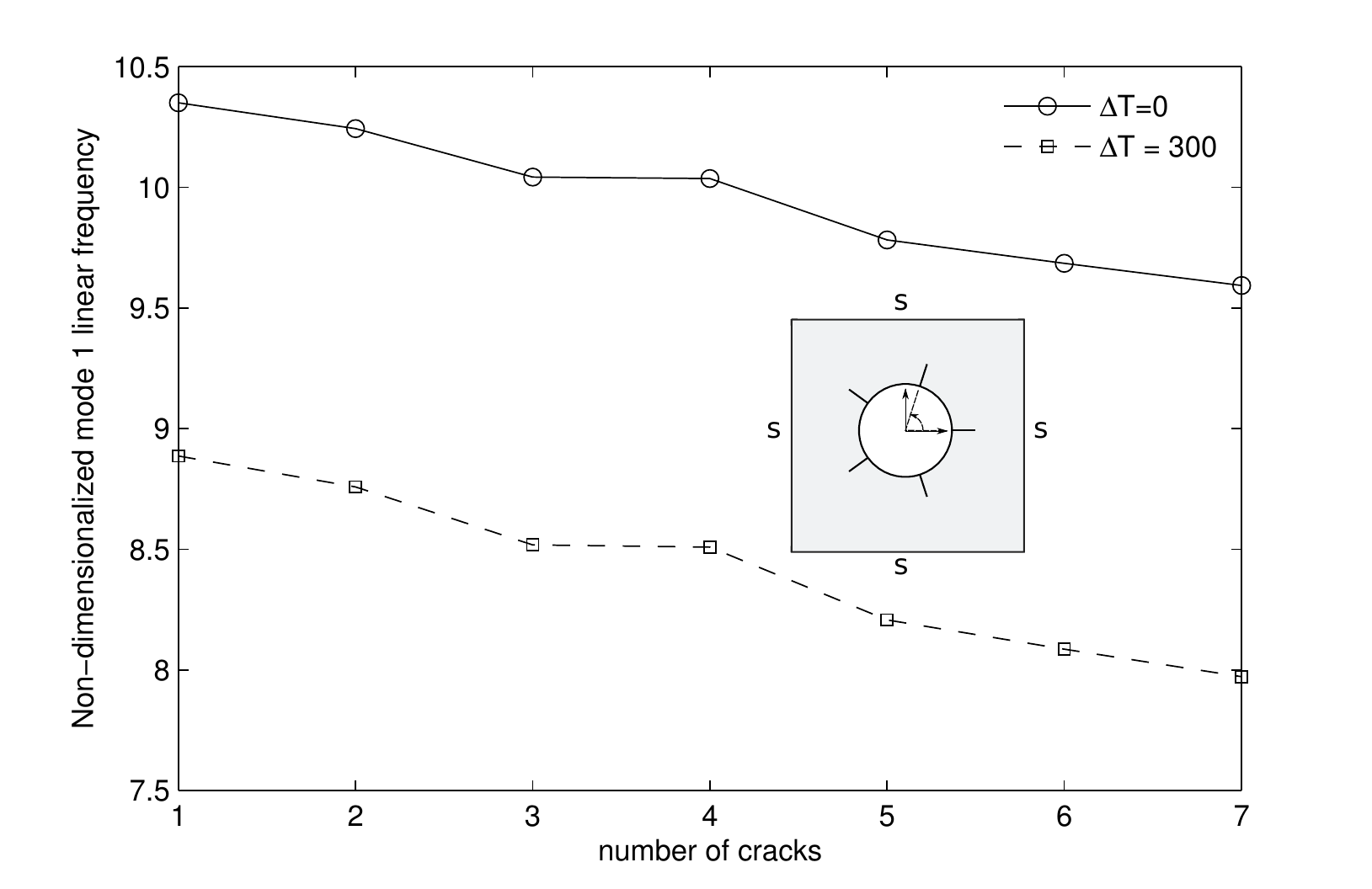}
\end{overpic}
\caption{The effect of the number of cracks on the mode 1 natural frequency $(\Omega)$ for a square simply supported FGM plate ($a/h=$ 10, $k=$2) with a circular cutout $r/a=0.2$ in different thermal environments. The length of each crack is $\ell/a=0.1$.}
\label{fig:NF_c1_7_n2_T600_Rc01_ssss}
\end{figure}

Next, the influence of multiple cracks emanating from a circular cutout is studied. \fref{fig:NF_c1c2_n2_T300_Rc01_ssss} shows the influence of the crack orientation on the first fundamental frequency. In this case, a simply supported square FGM plate with gradient index, $k=$ 2, a centrally located circular cutout $r/a=$ 0.2 and crack with length $\ell/a=$ 0.1 is considered. It can be seen from \fref{fig:NF_c1c2_n2_T300_Rc01_ssss} that the fundamental frequency is symmetric with respect to crack orientation $\theta=$ 45$^\circ$ and has a least value at $\theta=$ 45$^\circ$. \fref{fig:NF_c1_7_n2_T600_Rc01_ssss} shows the influence of number of cracks emanating from a circular cutout on the fundamental frequency for a simply supported square FGM plate with $a/h=$ 10, gradient index $k=$ 2, cutout ratio $r/a=$ 0.2 and exposed to different thermal gradients. It can be seen that with increasing number of cracks and thermal gradient, the fundamental frequency decreases. This can be attributed to the increase in the local flexibility due to the presence of cracks and thus decreasing the frequency.

\begin{figure}[ht!]
\centering
\begin{overpic}[scale=.75,unit=1mm]{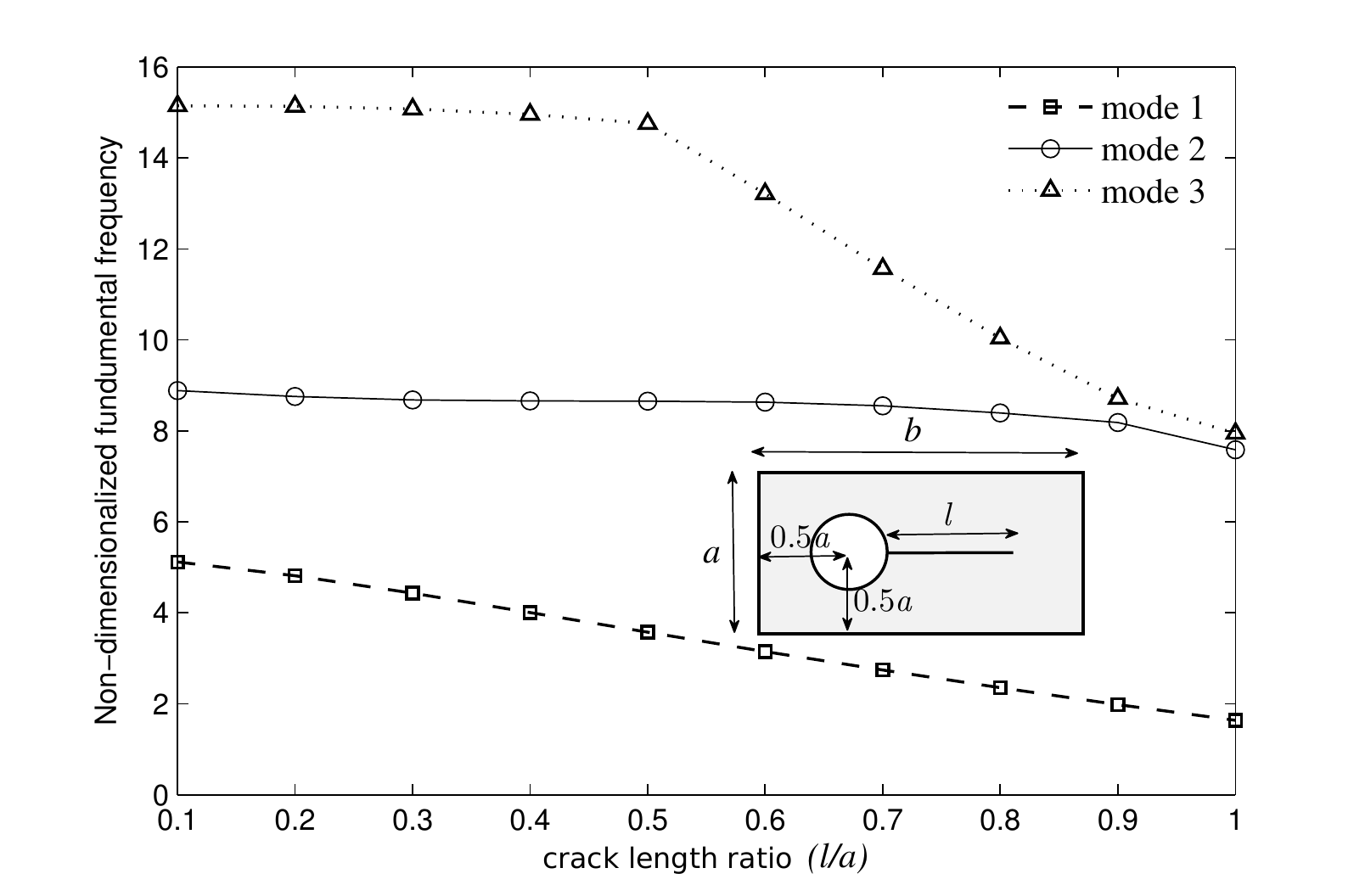}
\end{overpic}
\caption{The effect of the crack length $\ell/a$ on the first three fundamental frequencies $(\Omega)$ for a rectangular simply supported FGM plate ($b/a=$ 2, $a/h=$ 10, $k=$ 2) with a circular cutout $r/a=0.2$ in a thermal environment with $\Delta T=$300K.}
\label{fig:NF_c1_n2_T600_sss}
\end{figure}

\paragraph{Influence of the crack length} As a last example, the influence of the crack length on the first three fundamental frequency is studied. The geometry of the cutout, the location of the crack and the plate geometry is shown in \fref{fig:NF_c1_n2_T600_sss}. The plate is exposed to a thermal gradient, $\Delta T = $ 300K ($T_c = $ 600K and $T_m=$ 300K). It can be seen from \fref{fig:NF_c1_n2_T600_sss} that increasing the crack length decreases the fundamental frequencies. The decrease in the mode 1 frequency is linear whilst the mode 3 frequency decreases in a bi-linear mode. For the rectangular plate considered here, the effect of increasing crack length has minimal impact on the mode 2 frequency. This is because for the rectangular plate (see \fref{fig:NF_c1_n2_T600_sss}, the crack is perpendicular to the deformed shape.


%

\section{Summary}
The influence of a centrally located cutout and cracks emanating from the periphery of the cutout on the fundamental frequency of FGM plates is numerically studied. The influence of various parameters, viz., the plate thickness $h$, the plate aspect ratio $a/b$, the cutout geometry $d/e$, the cutout radius $r/a$, the crack length $\ell/a$ and the gradient index $k$ is studied for FGM plates in different thermal environment. The formulation is based on the first order shear deformation theory for plates and an enriched four-noded field consistent shear flexible element is used. The material is assumed to be temperature dependent and graded in the thickness direction. From a systematic parametric study, the following can be concluded:

\begin{itemize}
\item Increasing the cutout radius $r/a$ and the cutout geometry $d/e$ increases the frequency. This can be attributed to stiffness degradation.
\item Increasing gradient index $k$ decreases the natural frequency. This is due to the increase in the metallic volume fraction.
\item Increasing the slenderness ratio $a/h$ and the aspect ratio $a/b$ increases the frequency. 
\item Increasing the crack length $\ell/a$ and increasing the number of cracks, decreases the natural frequency.
\item Increasing the thermal gradient decreases the frequency.
\end{itemize}

It can be observed that the combined effect of increasing the cutout radius, the cutout geometry and decreasing gradient index is to increase in the frequency. In all the cases, the increase is due to stiffness degradation. In case of the gradient index, decreasing metallic volume fraction increases the stiffness and thus increasing the frequency. 

\section{ Acknowledgements} 
The authors would like to thank the financial support of the School of Engineering at Cardiff University for Ahmad Akbari Rahimabadi PhD. S Natarajan would like to acknowledge the financial support of the School of Civil and Environmental Engineering, The University of New South Wales for his research fellowship since Sep 2012. S Bordas would like to acknowledge the partial financial support of the Framework Programme 7 Initial Training Network Funding under grant number 289361 "Integrating Numerical Simulation and Geometric Design Technology". S Bordas also thank partial funding provided by: (1) the EPSRC under grant EP/G042705/1 Increased Reliability for Industrially Relevant Automatic Crack Growth Simulation with the eXtended Finite Element Method and (2) the European Research Council Starting Independent Research Grant (ERC Stg grant agreement No. 279578) entitled “Towards real time multiscale simulation of cutting in non-linear materials with applications to surgical simulation and computer guided surgery"

\newpage
\bibliographystyle{plain}
\bibliography{myRefFGM}

\end{document}